\documentclass{lms}
\usepackage{amssymb,amsmath,amstext,amsxtra,comment}
\usepackage{bm}       % allows bold italic letters in math mode
\usepackage{mathtools} % allows more extendible arrows
\usepackage{xcolor}
\usepackage[all]{xy}
\usepackage{booktabs}
\usepackage{colonequals}
\usepackage[breaklinks]{hyperref}
\hypersetup{colorlinks=true,urlcolor={blue!80!black},citecolor={blue!80!black},linkcolor={blue!80!black}}
\usepackage{enumerate}
\usepackage[group-minimum-digits=3,group-separator=\text{,}]{siunitx}
\usepackage{breakurl}

\numberwithin{equation}{subsection}

\newnumbered{assertion}[equation]{Assertion}    % 1st argument is your name for it
\newnumbered{conjecture}[equation]{Conjecture}  % 2nd argument is what is printed
\newnumbered{definition}[equation]{Definition}
\newnumbered{hypothesis}[equation]{Hypothesis}
\newnumbered{remark}[equation]{Remark}
\newnumbered{note}[equation]{Note}
\newnumbered{observation}[equation]{Observation}
\newnumbered{problem}[equation]{Problem}
\newnumbered{question}[equation]{Question}
\newnumbered{algorithm}[equation]{Algorithm}
\newnumbered{example}[equation]{Example}
\newunnumbered{notation}{Notation} % This is usually unnumbered
\newcommand{\Z}{\mathbb Z}
\newcommand{\Q}{\mathbb Q}
\newcommand{\R}{\mathbb R}
\newcommand{\F}{\mathbb F}

\newcommand{\A}{\mathbb A}
\newcommand{\C}{\mathbb C}
\newcommand{\PP}{\mathbb P}

\newcommand{\Gal}{\mathrm{Gal}}
\newcommand{\Aut}{\mathrm{Aut}}
\newcommand{\ST}{\mathrm{ST}}

\newcommand{\defi}[1]{\textsf{#1}} 	% for defined terms
\newcommand{\abs}[1]{\lvert{#1}\rvert} 	% for defined terms

\DeclareFontFamily{U}{wncy}{}
\DeclareFontShape{U}{wncy}{m}{n}{<->wncyr10}{}
\DeclareSymbolFont{mcy}{U}{wncy}{m}{n}
\DeclareMathSymbol{\Sha}{\mathord}{mcy}{"58}

\DeclareMathOperator{\M}{M}

\DeclareMathOperator{\GL}{GL}
\DeclareMathOperator{\USp}{USp}

\newcommand{\U}{\mathrm{U}}

\DeclareMathOperator{\SU}{SU}

\DeclareMathOperator{\ord}{ord}
\DeclareMathOperator{\repart}{Re}

\DeclareMathOperator{\Jac}{\mathrm{Jac}}
\DeclareMathOperator{\opchar}{\mathrm{char}}
\DeclareMathOperator{\disc}{\mathrm{disc}}

\DeclareMathOperator{\End}{\mathrm{End}}
\DeclareMathOperator{\Res}{\mathrm{Res}}

\newcommand{\frakp}{\mathfrak{p}}

\DeclareMathOperator{\Qbar}{\overline{\Q}}

\newcommand{\avs}[1]{{\color{cyan} \textsf{[[#1]]}}}

\newcommand{\jrs}[1]{{\color{blue} \textsf{[[#1]]}}}

\usepackage{hyperref}
\hypersetup{pdftitle={A database of genus 2 curves over the rational numbers},pdfauthor={Andrew R. Booker, Jeroen Sijsling, Andrew V.~Sutherland, John Voight, and Dan Yasaki}} 
\hypersetup{colorlinks=true,linkcolor=blue,anchorcolor=blue,citecolor=blue}

\title[A database of genus 2 curves]% end with percent
 {A database of genus 2 curves over the rational numbers} % This is the full title of the paper
\author[Booker, Sijsling, Sutherland, Voight and Yasaki]{Andrew R. Booker, Jeroen Sijsling, Andrew V. Sutherland, \\ John Voight and Dan Yasaki}

\classno{11G30, 14H45, 14Q05 (primary)}
\extraline{
The first author was supported by EPSRC grants
EP/H005188/1, EP/L001454/1 and EP/K034383/1.
The third author was supported by NSF grants DMS-1115455 and DMS-1522526. 
The fourth author was supported by an NSF CAREER Award (DMS-1151047).  
The fifth author was partially supported by NSA grant H98230-15-1-0228. This manuscript is submitted for publication with the understanding that the United States government is authorized to produce and distribute reprints.

}

\begin{document}
\maketitle

\begin{abstract}
% An \emph{abstract} written in English is
% required and should preferably have fewer than 200 words.
% Please do not include citations, footnotes or references
% to numbered equations, figures, tables or theorems in your
% abstract. Avoid complicated formulae or displayed equations, if
% possible.
We describe the construction of a database of genus $2$ curves of small discriminant that includes geometric and arithmetic invariants of each curve, its Jacobian, and the associated $L$-function.  This data has been incorporated into the $L$-Functions and Modular Forms Database (LMFDB).  
\end{abstract}

\section{Introduction}

The history of computing tables of elliptic curves over the rational numbers extends as far back as some of the earliest machine-aided computations in number theory.  The first of these tables appeared in the proceedings of the 1972 conference held in Antwerp~\cite{AntwerpIV}.  Vast tables of elliptic curves now exist, as computed by Cremona \cite{Cremona0,CremonaTables,Cr14} and Stein--Watkins \cite{SteinWatkins}.  These tables have been used extensively in research in arithmetic geometry to test and formulate conjectures, and have thereby motivated many important advances.

In this article, we continue this tradition by constructing a table of genus $2$ curves over the rational numbers.  We find a total of \num{66158} isomorphism classes of curves with absolute discriminant at most $10^6$; for each curve, we compute an array of geometric and arithmetic invariants of the curve and its Jacobian, as well as information about rational points and the associated $L$-function.  This data significantly extends previously existing tables (see section~\ref{sec:compl} for a comparison and discussion of completeness), and it is available online at the \emph{$L$-functions and Modular Forms Database (LMFDB)} \cite{lmfdb}.  We anticipate that the list of curves and the auxiliary data we have computed will be useful in ongoing investigations in genus $2$, including the Sato--Tate conjecture as formulated by Fit\'e--Kedlaya--Rotger--Sutherland \cite{FKRS12} and the paramodular conjecture of Brumer--Kramer \cite{BK14}.

This article is organized as follows.  In section \ref{sec:notation}, we set up basic notation and background; in section \ref{sec:boxsearch}, we describe our search procedure; in sections \ref{sec:invs}--\ref{sec:end}, we outline the computation of various invariants associated to the curve and its Jacobian; and 
% We conclude  with some observations and questions, and 
in the appendix, we provide tables exhibiting some interesting aspects of the data.
%\avs{Given the page breaking, we could easily add a few sentences to the abstract or introduction without impacting the length of the article if anyone has something else they would like to say here.}

\section{Notation and background} \label{sec:notation}

In this section, we set up basic notation and recall some relevant background material.
For a survey of computational aspects of genus 2 curves, see the article by Poonen \cite{Po94}.  
%\subsection{Notation} \label{sec:notation}

Let $A$ be an abelian variety over a field $F$.  For a field extension $K \supseteq F$, we denote by $A_K$ the abelian variety over $K$ obtained from $A$ by base extension from $F$ to $K$.
%  and we take number fields to be subfields of $\Qbar$.  
% For an abelian variety $A$ over $F$ and any algebraic extension $L$ of $k$ in $\Qbar$, we write $A_L\colonequals A\times_k L$ to denote the abelian variety over $L$ obtained from $A$ by base extension from $k$ to $L$.
% If $A$ and $B$ are abelian varieties over $k$, the notation $A\sim B$ indicates that $A$ and $B$ are related by an isogeny that is defined over $k$, in which case we say that $A$ and $B$ are isogenous; to indicate that $A$ and $B$ are isogenous over an algebraic extension $L/k$, we write $A_L\sim B_L$ and say that $A$ and $B$ are $L$-isogenous.
We write $\End(A)$ to denote the ring of endomorphisms of $A$ \emph{defined over $F$}; we write $\End(A_K)$ for endomorphisms of $A$ defined over an extension $K$.  For a $\Z$-algebra $R$ we use $\End(A)_R \colonequals  \End(A) \otimes_\Z R$ to denote the $R$-endomorphism algebra of $A$.  Finally, by a \defi{curve} over a field $F$ we mean a smooth, projective, geometrically integral scheme of finite type over $F$ of dimension $1$.  

\subsection{Models and discriminants}

Let $F$ be a perfect field and let $X$ be a curve of genus 2 over~$F$.  Then $X$ has a \defi{Weierstrass model} (or \defi{Weierstrass equation})
% It follows easily from the Riemann-Roch theorem\footnote{Fix an effective canonical divisor~$D$, pick any $x\in \mathcal{L}(D)-k$, any $y\in \mathcal{L}(3D)-{\rm span}\{1,x,x^2,x^3\}$, and note that $\mathcal{L}(6D)$ has dimension 11 and contains the 12 functions $1,x,\ldots,x^6,y,xy,\ldots,x^3y$.} 
\begin{equation} \label{eqn:weierstrass}
y^2+h(x)y=f(x)
\end{equation}
where $f,h \in F[x]$ have degrees bounded by $6$ and $3$, respectively.  It follows that $X$ is hyperelliptic over $F$, with the map $x:X \to \PP^1$ of degree $2$.
%  and the ring $F[x,y]/(y^2+h(x)y-f(x))$ is an integral domain whose fraction field is isomorphic to the function field $F(X)$ of $X$.  It follows that $X$ is \defi{hyperelliptic} (${\rm genus}(C)\ge 2$ and $[k(C):k(x)]=2$).

Now suppose that $\opchar F \ne 2$.  Completing the square yields a \defi{simplified} Weierstrass model $y^2=4f(x)+h(x)^2$ for $X$.
For $n\ge 1$, let $g_{n,\textup{univ}}(x)=a_nx^n + \dots + a_1x+a_0 \in \Z[a_0,\dots,a_n][x]$ and define
\[ \disc_n \colonequals  \frac{1}{a_n}(-1)^{n(n-1)/2} \Res(g_{n,\textup{univ}}(x),\frac{d}{dx}g_{n,\textup{univ}}(x)) \in \Z[a_0,\dots,a_n]. \]
For a polynomial $g \in F[x]$ with $\deg(g) \leq n$, define $\disc_n(g) \in F$ by specialization.  Following Liu \cite{Li94}, the \defi{discriminant} of the Weierstrass model $y^2+h(x)y = f(x)$ is defined by
\[
\Delta(f,h) \colonequals 2^{-12}\disc_6(4f+h^2).
\]
For $f_{\textup{univ}}(x)=\sum_{i=0}^6 f_i x^i$ and $h_{\textup{univ}}(x)=\sum_{i=0}^3 h_ix^i$ in $\Z[f_0,\ldots,f_6,h_0,\ldots,h_3][x]$, the greatest common divisor of the coefficients of $\disc_6(4f_{\textup{univ}}+h_{\textup{univ}}^2)$ is $2^{12}$, so 
\[
\Delta(f_{\textup{univ}},h_{\textup{univ}}) \in \Z[f_0,\dots,f_6,h_0,\dots,h_3],
\]
and in particular, if $f,h \in R[x]$, then $\Delta(f,h) \in R$ for a domain $R \subseteq F$.
For any $f,h\in F[x]$ with $\deg f\le 6$ and $\deg h\le 3$, we have $\Delta(f,h)\ne 0$ if and only if $y^2+h(x)y=f(x)$ is a Weierstrass model for a genus 2 curve over $F$.
In the special case that $f$ is monic with $\deg(f)=5$ and $\deg h \le 2$ (which can be achieved if and only if $X$ has a Weierstrass point defined over $F$), this agrees with the definition given by Lockhart \cite{Lo94}.
%In characteristic~2 one instead works with lifts of the polynomials $f,h$ to the ring of Witt vectors $W(k)$ and defines $\Delta(f,h)$ as above, modulo $2W(k)$.

If $(y')^2+h'(x')y'=f'(x')$ 
is another Weierstrass model for the curve $X$, then it is related to \eqref{eqn:weierstrass} via a change of variables of the form
\begin{equation}\label{eq:trans}
x' = \frac{ax+b}{cx+d},\quad y' =\frac{ey+j(x)}{(cx+d)^3}, \quad \text{with }  ad-bc,e\in F^\times,\ j\in F[x],\ \deg j \le 3,
\end{equation}
and we have $\Delta(f',h')=e^{20}(ad-bc)^{-30}\Delta(f,h)$; see Liu \cite[\S 1.3]{Li94}.
Thus, as an element of $F^\times/F^{\times 10}$, the discriminant $\Delta(f,h)$ is an invariant of $X$, and if $v$ is a discrete valuation on~$F$, then $v(\Delta(f,h))\bmod 10$ is also an invariant.

Now suppose that $F$ is a local field or a (global) number field, and let $R$ be the ring of integers of $F$.
An \defi{integral model} of $X$ is a Weierstrass model \eqref{eqn:weierstrass} with $f,h \in R[x]$, which implies $\Delta(f,h) \in R$.  We may rescale any given Weierstrass model to obtain an integral model after clearing denominators.
% We can always convert a given Weierstrass model to an integral model by multiplying through by $a^2$ for some $a\in R$ divisible by the denominator of every coefficient of $f$ or $h$ and then replacing $y$ by $y/a$.
If $F$ is a non-archimedean local field with discrete valuation $v$, then a \defi{minimal model} for $X$ is an integral model that minimizes $v(\Delta(f,h))\in \Z_{\ge 0}$, and the \defi{minimal discriminant} $\Delta_{\min}(X)$ is the $R$-ideal generated by $\Delta(f,h)$ for any minimal model $y^2+h(x)y=f(x)$.
If instead $F$ is a number field, then the minimal discriminant of $X$ is defined by
\[
\Delta_{\min}(X)\colonequals \prod_{\frakp} \Delta_{\min}(X_{F_{\frakp}}) \subseteq R,
\]
where $\frakp$ ranges over the non-zero prime ideals of $R$ and
$F_\frakp$ is the completion of~$F$ at $\frakp$, and we say that an
integral model for $X$ is a \defi{global minimal model} if
$\Delta(f,h)R=\Delta_{\min}(X)$.  Global minimal models do not
always exist, but if $F$ has class number one, then~$X$
has a global minimal model \cite[Remarque 6]{Li96}.  Minimal
models for genus 2 curves over~$\Q$ can be computed in \textsc{Magma}~\cite{magma} via the function \texttt{MinimalWeierstrassModel}, implemented by Michael Stoll.

When $F=\Q$, the discriminant $\Delta(f,h) \in \Z$ of every minimal model is the same, since it is determined up to sign by the principal ideal $\Delta_{\min}(X)$ and $-1\not\in\Q^{\times 10}$; we then view $\Delta_{\min}(X)$ as the integer $\Delta(f,h) \in \Z$ and define the \defi{absolute discriminant} of $X$ to be $\abs{\Delta_{\min}(X)} \in \Z_{>0}$.

\section{Searching in a box} \label{sec:boxsearch}

We wish to enumerate integral models $y^2+h(x)y=f(x)$ of genus 2 curves $X$ over $\Q$ with the goal of finding as many non-isomorphic curves as possible with absolute discriminant $\abs{\Delta_{\min}(X)}\le D$, where $D\colonequals 10^6$ is our chosen discriminant bound.
Computing the minimal discriminant of a given curve $X$ is time-consuming, so we instead just search for integral models with $\abs{\Delta(f,h)}\le D$.
We then minimize these models if necessary, but this is rarely required; absolute discriminants bounded by $10^6$ will necessarily have minimal valuation at all primes $p \geq 5$, since $5^{10}>10^6$.
Our approach is to enumerate a large finite set $S$ of integral models from which we will select those with $|\Delta(f,h)|\le D$.  The set $S$ is defined by constraints on the coefficients of $f,h\in \Z[x]$ that are chosen so that $S$ contains many minimal models.

\subsection{The shape of the box}

We can assume without loss of generality that the coefficients $h_i$ of $h(x)=\sum_{i=0}^3 h_i x^i$ lie in the set $\{0,1\}$: indeed, if $h'\equiv h\pmod 2$, applying \eqref{eq:trans} with $j=(h'-h)/2$ and $e=1$ and $x'=x$ yields an integral model $(y')^2+h'(x)y'=f'(x)$ with the same discriminant.
We cannot, however, assume $h=0$;
indeed, only a small fraction
of the curves we eventually found have minimal models with $h=0$.
This leaves $16$ possible $h\in \Z[x]$ with $\deg h\le 3$ to consider.
One can reduce this to 6 via the substitutions $x\mapsto 1/x$ and $x\mapsto x+1$ which leave $\Delta(f,h)$ fixed; however, we choose not to do this as these transformations may move a model outside the set~$S$ we are searching (we do use this symmetry to optimize our search).

For the coefficients of $f=\sum_0^6 f_ix^i$ we impose constraints of the form $|f_i|\le B_i$ that determine the shape of a box (orthogonal parallelepiped) in which we hope to find integral models with small discriminants. We experimented with various shapes and eventually settled on three:
\begin{align*}
S_1(B)&\colonequals \{(f,h):|f_i|\le B,\ h_i=0,1\},\\
S_2(a,b)&\colonequals \{(f,h):|f_i|\le ab^{6-i},\ h_i=0,1\},\\
S_3(b)&\colonequals \{(f,h):|f_i|\le b^{4-|i-3|},\ h_i=0,1\}.
\end{align*}
The ``flat'' box $S_1$ is a hypercube (a ball in the $L^\infty$-metric).  The ``weighted'' box $S_2$ complements~$S_1$ and tries to find curves with small discriminant but large ``flat'' height, with weighting designed to give each term in the discriminant roughly the same size.  The ``hill'' box $S_3$ complements $S_1,S_2$ and tries to find curves with small discriminant but large ``flat'' and ``weighted'' height: it allows larger middle coefficients with smaller end coefficients.  Taken together, the boxes $S_i$ are a crude attempt to take the oddly shaped region defined by $\disc(f,h) \leq X$ and approximate it with boxes that have little overlap.  Table~\ref{tab:boxsearch} lists the number of isomorphism classes with $|\Delta(f,h)|\le 10^n$ for $n=3,4,5,6$, using boxes of cardinality $10^m$ with $m=11,\ldots,17$ for these three shapes.  

To squeeze out a bit more, we also considered sets
\[
S_4(b,d)\colonequals \{(f,h): \textstyle{\sum}_i \lceil\log_b(|f_i|+1)\rceil\le d,\ h_i= 0,1\},
\]
which can be viewed as a union of $\binom{d+6}{6}$ boxes.  For our final search we settled on 
\[
S_1 \colonequals S_1(90),\qquad S_2 \colonequals S_2(2,3.51),\qquad S_3 \colonequals S_3(7.14),\qquad S_4 \colonequals S_4(10,10),
\]
with cardinalities $\#S_1,\#S_2,\#S_3\approx 10^{17}$ and $S_4\approx 10^{16.3}$; see Table~\ref{tab:boxsearchfinal} for more details.

\subsection{Computing discriminants}\label{sec:compdisc}

The discriminant polynomial $\disc_6\in \Z[a_0,\ldots,a_6]$ is a homogeneous polynomial of degree~10 with 246 terms.
Directly computing $\Delta(f,h)=2^{-12}\disc_6(4f+h^2)$ for a single model $(f,h)$ involves at least 1000 arithmetic operations, but to compute $\Delta(f,h)$ for all $(f,h)$ in a box we instead use an algorithm that efficiently enumerates values of a multivariate polynomial $g(x_1,\ldots,x_n)$ over integer points in a box $A_1\times\cdots\times A_n$ of $\Z^n$.
With this approach we require less than 6 clock cycles, on average, to compute each discriminant $\Delta(f,h)$ for $(f,h)\in S$.

The algorithm first creates a \defi{monomial tree} with $n$ levels.
A node at level $m\le n$ in the tree corresponds to a monomial in $x_1,\ldots,x_m$, identified by a vector of $m$ non-negative exponents, and contains storage for an integer coefficient.
The nodes at level $n$ (the leaves) correspond to the terms of $g$.
The nodes at levels $m < n$ correspond to the terms of a partially instantiated polynomial $g(x_1,\ldots,x_m,a_{m+1},\ldots,a_n)$, for some $a_{m+1},\ldots,a_n\in \Z$; the monomial exponents depend only on $g$, but the coefficients vary with $a_{m+1},\ldots,a_n$.
Each node at level $m+1$ has an edge leading to the node at level $m$ corresponding to the substitution $x_{m+1}=a_{m+1}$.

For $n=1$, the polynomial $g(x_1)$ is univariate, and we may efficiently enumerate its values at integers in the interval $A_1$ using finite differences, as in Kedlaya--Sutherland~\cite[\S 3]{KS08}.  Each step requires just $d$ additions, where $d=\deg_{x_1} g$.
For $n>1$ we reduce to the case $n-1$ as follows.
We begin with $a_n=\min (A_n)$ and compute the powers $a_n^i$ for $0\le i \le d$, where $d=\deg_{x_n} g$.
We then use these powers to instantiate the monomials at level $n$ of the tree at $x_n=a_n$, using one multiplication per node, thereby obtaining a coefficient value for each instantiated monomial.
We then set the coefficient value of each node at level $n-1$ to the sum of the coefficient values of its instantiated children.
The nodes at level $n-1$ now represent the terms of the polynomial $g(x_1,\ldots,x_{n-1},a_n)$ in $n-1$ variables whose values we may recursively enumerate over integer points in $A_1\times\cdots\times A_{n-1}$.
We then increment $a_n$ and repeat until $a_n> \max (A_n)$.
%The sets $A_2,\ldots, A_n$ need not be intervals, but we assume $A_1$ is.

The time to enumerate the values of $g(x_1,\ldots,x_n)$ over $A_1\times \cdots\times A_n$ is dominated by the time spent in the bottom layer of the recursion ($n=1$), which involves $\deg_{x_1} g$ integer additions.
We treat each of the 16 possible values of $h$ as fixed, and for each $h$ view $\disc_6(4f_{\rm univ}+h^2)$ as a polynomial in the $n=7$ variables $f_0,\ldots,f_6$.
Using a monomial tree to enumerate values of $g \colonequals \disc_6(4f_{\rm univ}+h^2)\in \Z[f_0,\ldots,f_6]$ yields an algorithm that essentially uses just $\deg_{f_0} g = 5$ additions per discriminant.  As a further optimization, we work modulo $2^{64}$ so that each addition can be accomplished in a single clock cycle.
The monomial tree for $g$ contains 703 nodes and requires less than 8 kilobytes of storage, which easily fits in the L1 data cache of a modern~CPU.

To find integral models $y^2+h(x)y=f(x)$ with small discriminant, 
we look for pairs $(f,h)$ for which $\Delta(f,h)\equiv \Delta \pmod{2^{64}}$ for some integer $\abs{\Delta}\le D$.
For the relatively few pairs that pass this test, we recompute $\Delta(f,h)$ over $\Z$ to check if in fact $\abs{\Delta(f,h)} \leq D$;
this recomputation occurs so rarely that it has a negligible impact on the total running time.

\subsection{Parallel computation}\label{sec:comp}

Our search can be trivially parallelized by partitioning each box into products of the form $A_1\times\cdots\times A_n$ and applying the algorithm described in the previous section.
Using the compute engine on Google's Cloud Platform \cite{Google}, we ran a parallel computation to enumerate all pairs $(f,h)$ with $|\Delta(f,h)|\le 10^6$ in the sets $S_1,S_2,S_3,S_4$ defined above, which together contain over $3\times 10^{17}$ pairs $(f,h)$.
We allocated \num{2250} high-CPU compute nodes distributed across six of Google's compute engine regions (two in the central U.S., three in the eastern U.S., and one in western Europe).
Each high-CPU node comes equipped with 32 Intel Haswell cores and 28.8~GB of memory, yielding a total of \num{72000} cores.
We performed approximately 50 core-years of computation over the course of less than 8 hours of elapsed time.
To reduce the cost we used preemptible nodes.
On average, roughly 3 percent of our nodes were preempted each hour; these were automatically restarted from checkpoints that were created every few minutes.

\subsection{Completeness} \label{sec:compl}

We now compare our results with existing tables of genus 2 curves over $\Q$.  In every case we find that our tables include all previously known examples that fall within our discriminant range (as well as many that are completely new).

\begin{itemize}
\item The most comprehensive list of genus 2 curves over $\Q$ previously compiled is due to Stoll~\cite{St13}, who searched for curves with odd absolute discriminants bounded by~$11^4$, finding \num{823} distinct $\Q$-isomorphism classes with \num{749} distinct discriminants.
Our search found all of these curves, as well as many more: in total, we found \num{988} $\Q$-isomorphism classes of genus~2 curves over $\Q$ with \num{852} distinct odd discriminants of absolute value up to $11^4$.
Our tables differ already at the third smallest absolute minimal discriminant $|\Delta_{\min}(X)|=277$, where we find two non-isomorphic curves
\begin{align*}
&\href{http://www.lmfdb.org/Genus2Curve/Q/277/a/277/1}{y^2+(x^3+x^2+x+1)y=-x^2-x},\\% hack to avoid & in href -- AVS
&\href{http://www.lmfdb.org/Genus2Curve/Q/277/a/277/2}{y^2+y=x^5 - 9x^4 + 14x^3 - 19x^2 + 11x - 6},
\end{align*}
only one of which (the first) was found by Stoll \cite{St13}.
The smallest odd absolute discriminant found in our search that is new is \num{1665}, which arises for the curve
\[
\href{http://www.lmfdb.org/Genus2Curve/Q/1665/a/1665/1}{y^2+(x^3+x^2+1)y = x^4 + x^3 + 2x^2 + 2x + 1}.% this model is not the same as the one currently in the LMFDB, but that will change -- AVS
\]

\item Merriman--Smart \cite{MS93,Sm97} determined all genus 2 curves over $\Q$ with good reduction away from $2$.
Up to $\Q$-isomorphism there are 427 such curves, of which \href{http://www.lmfdb.org/Genus2Curve/Q/?abs_disc=128,256,512,1024,4096,8192,16384,32768,65536,131072,262144,524288}{29 curves} have absolute discriminant bounded by $10^6$,
all found by our search.  (The models given by Smart \cite[Table]{Sm97} are not necessarily minimal; we used Magma \cite{magma} to compute minimal models.)

\item Baker--Gonz\'alez--Gonz\'alez-Jim\'enez--Poonen \cite{BGGP05} and Gonz\'alez--Gonz\'alez-Jim\'enez \cite{GG02} determined a complete list of new modular curves of genus $2$, i.e., 
% JV: latin abbreviations (etc. e.g.) are not usually italicized
those curves $X$ that admit a dominant morphism $\pi\colon X_1(N)\to X$ from the modular curve $X_1(N)$ 
(for which $\pi^*(\Jac(X))\subseteq J_1(N)_{\textup{new}}$) for some level $N\in \Z_{\geq 1}$.  
Up to $\Q$-isomorphism there are~213 such curves: 149 have $\Q$-simple Jacobians \cite{GG02} and 64 have split Jacobians \cite{BGGP05}.
Of these 213 curves, all 61 with absolute discriminant bounded by $10^6$ were found in our search.

\item Gon\-z\'alez--Gonz\'alez-Jim\'enez\---Gu\`ardia \cite{GGG02} (see also Gonz\'alez--Gu\`ardia--Rotger \cite{GGR05}) computed equations of genus $2$ curves whose Jacobians are modular abelian varieties, i.e., they are subabelian varieties of the Jacobian $J_1(N)$ for some $N \in \Z_{\geq 1}$.  Among the 75 curves they list, the 3 with absolute discriminant bounded by $10^6$ were found in our search. 

\item Gonz\'alez--Rotger \cite{GR04} computed explicit models for some Shimura curves of genus 2; of these, the $6$ curves with absolute discriminant bounded by $10^6$ were also found in our search.

\item In pursuit of the paramodular conjecture, Brumer--Kramer \cite{BK14} list 38 genus 2 curves $X$ whose Jacobians $\Jac(X)$ have $\End(\Jac(X))=\Z$ of odd conductor less than $1000$.  Of these, $33$ have absolute discriminant bounded by $10^6$, all of which were found in our search.
In total we found 56 $\Q$-isomorphism classes of genus 2 curves $X$ whose Jacobians have odd conductor less than $1000$, of which \href{http://www.lmfdb.org/Genus2Curve/Q/?cond=169,249,277,295,349,353,363,389,427,461,523,529,555,587,597,603,691,709,713,731,741,743,745,763,797,807,841,847,893,909,925,953,961,971,975,997&disc=&num_rat_wpts=&torsion=&torsion_order=&two_selmer_rank=&is_gl2_type=False&st_group=&real_geom_end_alg=&aut_grp=&geom_aut_grp=}{43 curves} have $\End(\Jac(C))=\Z$.
These all appear to lie in the $33$ isogeny classes identified by Brumer--Kramer using Jacobians of curves with absolute discriminants below $10^6$.

\item In a similar vein, Farmer--Koutsoliotas--Lemurell \cite{fkl} used analytic methods to determine a complete list of the integers $N\le 500$ that can arise as the conductor of degree-4 $L$-functions under a standard set of hypotheses.
Excluding $L$-functions arising from products of elliptic curves or classical modular forms with real quadratic character, they list 12 values of $N\le 500$ that may arise for the $L$-function of a genus 2 curve \cite[Theorem 2.1]{fkl} (the 9 odd values match the Brumer--Kramer list up to $500$).
We found \href{http://www.lmfdb.org/Genus2Curve/Q/?cond=249,277,295,349,353,388,389,394,427,461,464,472}{21 curves} matching these 12 values (and at least one for each $N$), including \href{http://www.lmfdb.org/Genus2Curve/Q/?cond=388,464,472}{6 curves} that match the even values $388$, $464$, $472$ on their list.
\end{itemize}

Despite our attempts at exhaustiveness, the data in Table~\ref{tab:boxsearch} suggest that our tables are still incomplete.
Indeed, there are genus~2 curves with absolute minimal discriminant bounded by $10^6$ that do not have an integral model in any of the sets $S_1,S_2,S_3,S_4$ that we searched: we found 28 such curves while experimenting with boxes of different shapes.
One way to get a sense for the completeness of our tables is to consider an analogous search among genus 1 curves.
The discriminant polynomial then has degree 6, rather than 10, so we use the discriminant bound $\lfloor (10^6)^{6/10} \rfloor=3981$.
The Jacobian of such a genus~1 curve is isomorphic to an elliptic curve whose conductor must also be bounded by $3981$ and can therefore be found in Cremona's tables \cite{Cr14} (by Ogg's formula, the conductor of an elliptic curve over $\Q$ divides its minimal discriminant).
In total there are 1154 $\Q$-isomorphism classes of elliptic curves with absolute discriminant bounded by $3981$.
A search for genus 1 curves $y^2+h(x)y=f(x)$ with $\deg h \le 2$ and $\deg f\le 4$ using the same coefficient bounds as in the box $S_1(90)$ already hits 1132 of these $\Q$-isomorphism classes, or about 98 percent.
% \avs{Strictly speaking we are counting Jacobians here, not curves, but I think this is clear.}  % Yes, I think this is clear.

\section{Computing basic data}\label{sec:invs}

\subsection{Invariants}

The isomorphism class of a genus $2$ curve over an algebraically closed field can be characterized by a number of \emph{invariants}.  For elliptic curves, only the elliptic $j$-invariant is required.  In genus $2$, there are different choices of invariants available, all with their own distinct advantages and drawbacks.  We refer to seminal work by Igusa \cite{Igusa} and later work by Mestre~\cite{Mestre} for some theoretical expositions. The details of the implementation in \textsc{Magma} can be found in its Handbook \cite{MagmaHandbook}; the package was written by Everett Howe, based on routines developed by Fernando Rodriguez-Villegas, and with additional functionality added by Reynald Lercier and Christophe Ritzenthaler.

Let $F$ be a field with $\opchar F \neq 2$ and let $X$ be a genus $2$ curve over $F$.  Then $X$ has a simplified Weierstrass model
\[
  X : y^2 = f(x) = c \prod_{i = 1}^6 (x - \alpha_i) ,
\]
% \jv{To be completely in agreement with previous notation, we should write this as $y^2=g(x)$; but I prefer to think ``we may take $h(x)=0$'', so I hope this is OK}
where the $\alpha_i$ are the distinct roots of $f$ in an algebraic closure $\overline{F}$ of the base field $F$. Given two such roots $\alpha_i, \alpha_j$, we denote the difference $\alpha_i - \alpha_j$ by $(i j)$. We then form the expressions
\begin{equation}
  \begin{aligned}
    I_2 & \colonequals (4c)^2 \sum (1 2)^2 (3 4)^2 (5 6)^2 , \\
    I_4 & \colonequals  (4c)^4 \sum (1 2)^2 (2 3)^2 (3 1)^2 (4 5)^2 (5 6)^2 (6 4)^2 , \\
    I_6 & \colonequals  (4c)^6 \sum (1 2)^2 (2 3)^2 (3 1)^2 (4 5)^2 (5 6)^2 (6 4)^2 (1 4)^2
      (2 5)^2 (3 6)^2 , \\
    I_{10} & \colonequals  (4c)^{10} \prod (1 2)^2 = \disc_6(4f),
  \end{aligned}
\end{equation}
where each sum and product runs over the distinct expressions obtained by permuting the index set $\{1, \dots ,6\}$. The invariants $I_2, I_4, I_6, I_{10}$, defined by Igusa \cite[p.~620]{Igusa} by modifying a set of invariants due to Clebsch, are known as the \defi{Igusa-Clebsch invariants}. Two curves $X$ and $X'$ over $F$ are isomorphic over the algebraic closure $\overline{F}$ if and only if \[ (I_2:I_4:I_6:I_{10})=(I_2':I_4':I_6':I_{10}') \in \PP(2,4,6,10)(\overline{F}) \] describe the same point in weighted projective $(2, 4, 6, 10)$-space; that is, if and only if there exists a $\lambda \in \overline{F}^\times$ such that $(I_2',I_4',I_6',I_{10}') = (\lambda^2 I_2, \lambda^4 I_4, \lambda^6 I_6, \lambda^{10} I_{10})$.  

The Igusa--Clebsch invariants break down in characteristic $2$, 
% where hyperelliptic curves instead have to be described in terms of Artin--Schreier covers. 
and to deal with this problem, Igusa \cite[pp.~617ff]{Igusa} found a common normal form for hyperelliptic curves in arbitrary characteristic. Motivated by this, he defined the invariants \cite[pp.~621--622]{Igusa}
\begin{equation}
  \begin{aligned}
    J_2    & \colonequals I_2 / 8, \\
    J_4    & \colonequals (4 J_2^2 - I_4) / 96, \\
    J_6    & \colonequals (8 J_2^3 - 160 J_2 J_4 - I_6) / 576, \\
    J_8    & \colonequals (J_2 J_6 - J_4^2) / 4, \\
    J_{10} & \colonequals I_{10} / 4096 , 
  \end{aligned}
\end{equation}
now called the \defi{Igusa invariants}, with $(J_2:J_4:J_6:J_8:J_{10}) \in \PP(2,4,6,8,10)(\overline{F})$.  Moreover, these invariants are compatible with specialization of the model, for example, reducing modulo a prime. 

While the Igusa invariants are in a sense the best available with respect to arithmetic questions, we need $5$ of them to describe a moduli space of dimension $3$. 
% \jv{Don't we want to define the absolute Igusa(--Clebsch) invariants?  I think this is what the cryptographers predominantly use.} \jrs{Yes, but this coincides with what is written below, with the exception that G2 is a bit more careful; cryptographers only typically only use the first formula.} 
However, the moduli space of genus $2$ curves over a given field admits a simpler description: the corresponding moduli space can always be described by $3$ affine coordinates. Such coordinates were given by Cardona--Quer--Nart--Pujolas \cite{CQNP2,CQNP1}.  For $\opchar F \neq 2$, these \defi{G2-invariants} (sometimes also known as \defi{absolute Igusa invariants}) can be defined by
 \begin{equation}
  (g_1, g_2, g_3) =
  \begin{cases}
    (J_2^5 / J_{10}, J_2^3 J_4/ J_{10}, J_2^2 J_6 / J_{10}), & \text{if $J_2 \neq 0$;} \\
    (0, J_4^5 / J_{10}^2, J_4 J_6 / J_{10}), & \text{if $J_2 = 0, J_4 \neq 0$;}\\
    (0, 0, J_6^5 / J_{10}^3), & \text{otherwise.}
  \end{cases}
\end{equation}
The G2-invariants are \emph{absolute}, with $(g_1,g_2,g_3) \in \A^3(F)$ in affine space.   This time, $X$ has potentially good reduction at a prime if and only if none of the G2-invariants reduce to $\infty$ modulo that prime.  Conversely, given three G2-invariants in a field $k$, there always exists a curve with these invariants. However, writing down such a curve may require a quadratic extension of $k$. This obstruction was studied by Mestre \cite{Mestre} and Cardona--Quer \cite{CQNP1}.

\subsection{$L$-factors at good primes}

The $L$-series of a genus $2$ curve $X$ over $\Q$ is a Dirichlet series with an Euler product
\[ L(X,s) = \sum_{n=1}^{\infty} \frac{a_n}{n^s} = \prod_p L_p(p^{-s}) \]
such that for primes $p$ of good reduction we have 
\[ Z(X/\F_p,T)=\exp\left(\sum_{r=1}^{\infty} \#X(\F_{p^r})\frac{T^r}{r}\right)=\frac{L_p(T)}{(1-T)(1-pT)} \] 
and $L_p(T) \in 1+T\Z[T]$ has degree $4$.
Efficient methods for computing the coefficients of $L_p(T)$ for all good primes $p$ up to a given bound are discussed in Kedlaya--Sutherland \cite{KS08} and Harvey--Sutherland \cite{HS14,HS16} and available in \texttt{smalljac} \cite{smalljac}.

\subsection{Hash function}

If abelian varieties $A,A'$ over $\Q$ are isogenous over $\Q$, then $a_p(A)=a_p(A')$ 
for all primes $p$.
To facilitate the efficient grouping of curves into candidate isogeny classes (of their Jacobians), we define the \defi{hash} of an abelian surface $A$ to be the unique integer $h(A)\in [0,P-1]$ for which
\[
h(A)\equiv \sum_{\substack{2^{12} < p < 2^{13} \\ \text{$p$ prime}}} c_p a_p(A) \bmod P,
\]
where $P=2^{61}-1$ and $c_p=\lfloor \pi P^{e_p} \rfloor$, with $e_p \colonequals \#\{q : 2^{12}< q \leq p, \text{ $q$ prime}\}$ and $\pi=3.14\ldots$ (so the $c_p$ correspond to successive digits in the base-$P$ expansion of $\pi$).
Isogenous abelian surfaces necessarily have the same hash, and while the converse need not hold in general, it almost always does.
Hashes also provide a way to provisionally identify automorphic objects (such as modular forms) that have the same $L$-function.

\subsection{Sato--Tate group} \label{sec:st_group}

For an abelian surface $A$ over a number field~$F$, there are 52 distinct Sato--Tate groups $\ST(A)$ that arise, of which 34 occur when $F=\Q$ \cite[Theorem\ 1.5]{FKRS12}.  A coarser classification is obtained if we restrict our attention to the identity component $\ST(A^0)$, equivalently, the Sato--Tate group $\ST({A_K})$ of the base extension of $A$ to the minimal field $K$ over which all its endomorphisms are defined.
In this case there are just 6 possibilities, corresponding to 6 types of real endomorphism algebras $\End(A_K)_\R \colonequals \End(A_K)\otimes_\Z\R$, all of which can arise when $F=\Q$.  These are listed in Table \ref{table:connectedST}, along with the possible types of abelian varieties $A_K$ (up to isogeny).

\begin{table}[ht!] 
\begin{tabular}{llll} 
type & $\ST(A)^0$ & $\End(A_K)_\R$ & $A_K$\vspace{2pt}\\\toprule
$\mathbf{A}$&$\USp(4)$ & $\R$ & generic abelian surface\vspace{4pt}\\
$\mathbf{B}$&$\SU(2)\times \SU(2)$ & $\R\times\R$ & RM abelian surface or\\
&&& product of non-isogenous generic elliptic curves\vspace{4pt}\\
$\mathbf{C}$&$\U(1)\times \SU(2)$ & $\C\times \R$ & product of generic elliptic curve and CM elliptic curve\vspace{4pt}\\
$\mathbf{D}$&$\U(1)\times \U(1)$ & $\C\times\C$ & CM abelian surface or\\
&&&product of non-isogenous CM elliptic curves\vspace{4pt}\\
$\mathbf{E}$&$\SU(2)$ & $\mathrm{M}_2(\R)$ & QM abelian surface or\\
&&& square of a generic elliptic curve\vspace{4pt}\\
$\mathbf{F}$&$\U(1)$ & $\mathrm{M}_2(\C)$ & square of an elliptic curve with CM\vspace{2pt}\\
\bottomrule
\end{tabular}
\caption{Sato--Tate identity components and real endomorphism algebras of abelian surfaces.}\label{table:connectedST}
\end{table}

Under the generalized Sato--Tate conjecture, the Sato--Tate group $\ST(A)$ can be identified heuristically by computing the moments of~$a_p$ and $a_{p^2}$ to sufficiently high precision and consulting the corresponding tables \cite[Tables 9--10]{FKRS12}.
To obtain a rigorous identification that does not depend on the Sato--Tate conjecture, one computes the Galois module $\End(A_K)_\R$ (with the natural action of $\Gal(K/F)$) and applies the explicit correspondence with Sato-Tate groups \cite[Theorem~4.3]{FKRS12}.

\section{Conductor and bad Euler-factors}\label{sec:cond}

In this section, we describe the computation of the conductor $N\in \Z$ of an abelian surface $A/\Q$ that is the Jacobian of a genus 2 curve $X/\Q$ with minimal discriminant $\Delta_{\min}$; for further background and the definition of the conductor, see Brumer--Kramer \cite{BK94}.

\subsection{Computing the conductor}

Only primes dividing $\Delta_{\min}$ can divide the conductor $N$ (but the converse need not hold), as shown by Liu \cite{Li94}.  For an \emph{odd} prime $p \mid \Delta_{\min}$, the exponent $\ord_p(N)$ and the $L$-factor $L_p(T)$ can be computed using the algorithm of Liu \cite{Li94-alg} computing the stable reduction of the curve, implemented in Pari by Liu and Henri Cohen and in \textsc{Magma} by Tim and Vladimir Dokchitser.  There is also further recent work of Bouw--Wewers \cite{BouwWewers} to compute these invariants using the semistable model instead, and in theory this method will also compute $\ord_2(N)$ as well (the ``wild'' case).  However, a general implementation of this method is not yet available.

\subsection{Analytically computing the conductor}

So it remains to compute $\ord_2(N)$, and we do so analytically, using
methods that go back at least to Dokchitser \cite{D04} and the first
author \cite{B05}; see also recent work of Farmer--Koutsoliotas--Lemurell \cite{fkl}. 
The result of our computation will be that \emph{if}
the genus $2$ curve has an $L$-function with analytic continuation and
functional equation, then there is \emph{provably} only one possibility
for the conductor and bad $L$-factor.
One could take this unique possibility as the
definition of $\ord_2(N)$ and the $L$-factor, and it must agree with
the standard definition if the Hasse--Weil conjecture is true.
Moreover, by multiplicity $1$ for the Selberg class \cite{So04},
there is never more than one permissible choice of the bad $L$-factor, so
a suitably general version of the method will always succeed.

Let $A$ be an abelian surface over $\Q$ and let 
\[ L(A,s)=\prod_p L_p(A,p^{-s})^{-1}=\sum_{n=1}^\infty a_nn^{-s} \]
be its Hasse--Weil $L$-function.  Let
$\Lambda(A,s)=\Gamma_\C(s)^2L(A,s)$, where
$\Gamma_\C(s)=2(2\pi)^{-s}\Gamma(s)$.
Then conjecturally $\Lambda(A,s)$ continues to an entire function of
order $1$ and satisfies the functional equation
\begin{equation}
\label{LJfunceq}
\Lambda(A,s)=wN^{1-s}\Lambda(A,2-s),
\end{equation}
where $N \in \Z_{\geq 1}$ is the conductor of $A$ and $w=\pm1$ is its root number.

For $x>0$ we define
$$
S(x)=\frac1x\sum_{n=1}^\infty a_nK_0(4\pi\sqrt{n/x}),
$$
where $K_0(y)=\int_0^\infty e^{-y\cosh{t}}\,\textup{d}t$ is the $K$-Bessel
function. Recall that $K_0(y)\sim\sqrt{\frac{\pi}{2y}}e^{-y}$
as $y\to\infty$ and $K_0(y)\sim\log\frac1y$ as $y\to0^+$; thus,
$S(x)$ decays rapidly as $x\to0^+$, and from the
estimate $a_n\ll n^{\frac12+\varepsilon}$ we have 
$S(x)\ll x^{\frac12+\varepsilon}$ for large $x$. Therefore, the integral
$$
I(s)=\int_0^\infty S(x)x^{-s}\,\textup{d}x
$$
is absolutely convergent for $\Re(s)>\frac32$.
For $s$ in that region we may also change the order of sum and integral
and make the substitution $x\mapsto(4\pi/y)^2n$ to obtain
\begin{align*}
I(s)&=\sum_{n=1}^\infty a_n\int_0^\infty
K_0(4\pi\sqrt{n/x})x^{-s}\,\frac{\textup{d}x}x
=\sum_{n=1}^\infty a_n
\int_0^\infty 2\bigl((4\pi/y)^2n\bigr)^{-s}K_0(y)\,\frac{\textup{d}y}y\\
&=2(4\pi)^{-2s}L(A,s)\int_0^\infty K_0(y)y^{2s-1}\,\textup{d}y.
\end{align*}
Finally, from the Mellin transform identity
$\int_0^\infty K_0(y)y^{s-1}\,\textup{d}y=2^{s-2}\Gamma(s/2)^2$, we arrive at
$$
I(s)=\tfrac18\Lambda(A,s)
\quad\text{for }\repart(s)>\tfrac32.
$$

By a classical argument of Hecke, it then follows 
that the conjectured analytic
continuation, growth properties and functional equation
\eqref{LJfunceq} of $\Lambda(A,s)$ are equivalent to
the identity
\begin{equation}\label{Ssymmetry}
S(x)=wS(N/x)\quad\text{for all }x>0.
\end{equation}
Given a purported sequence of Dirichlet coefficients $\{a_n\}$ and
values for $N$ and $w$, this gives a falsifiable identity that we can
use to test for consistency with the Hasse--Weil conjecture.
To carry this out in practice, we fix a suitable constant $C>0$ and
consider the truncated sum
$$
S_C(x)=\frac1x\sum_{n\le Cx}a_nK_0(4\pi\sqrt{n/x}).
$$
Using the Ramanujan bound for $a_n$, one obtains the following estimate
for the truncation error, provided that $C\ge5$:
\begin{equation}\label{SCerror}
|S(x)-S_C(x)|<4Cx^{\frac14}\big(1+2x\sqrt{C}\big)e^{-4\pi\sqrt{C}}.
\end{equation}
Also, since we might not know the values of $a_n$ for even $n$,
we use the fact that $L(A,s)$ is given
by an Euler product to write
\begin{equation}\label{SCodd}
S_C(x)=\frac1x\sum_{j=0}^\infty a_{2^j}
\sum_{\substack{m\le Cx2^{-j}\\ 2\nmid m}}a_m K_0(4\pi\sqrt{m2^j/x})=
\sum_{j=0}^{\lfloor\log_2(Cx)\rfloor}
\frac{a_{2^j}}{2^j}S_C^{\rm odd}\!\left(\frac{x}{2^j}\right),
\end{equation}
where 
$$
S_C^{\rm odd}(x)=\frac1x\sum_{\substack{n\le Cx\\2\nmid n}}
a_nK_0(4\pi\sqrt{n/x}).
$$

We then simply try all possibilities for $N$, $w$, and the Euler
factor $L_2(A,T)$ to see which, if any, give answers consistent with
\eqref{Ssymmetry}. Specifically, we compute
$S(2^{\frac14}\sqrt{N})-wS(2^{-\frac14}\sqrt{N})$
for each candidate set of parameters and check whether any of
these are zero.  By Brumer--Kramer \cite[Theorem~6.2]{BK94} we have $\ord_2(N)\le20$.
Moreover, in the case that $A$ is the Jacobian of a genus $2$ curve
of discriminant $\Delta_{\min}$, by Liu \cite{Li94} we have $N\mid\Delta_{\min}$,
so we may assume that $\ord_2(N)\le\min(20,\ord_2(\Delta_{\min}))$.  Thus,
by \eqref{SCerror} and \eqref{SCodd}, our computation reduces to that
of $S_C^{\rm odd}(\sqrt{2^{k-\frac12}N_{\rm odd}})$ for each $k$ with
$0\le k\le\min(20,\ord_2(\Delta_{\min}))+1$, where $N_{\rm odd}$ denotes the
odd part of the conductor.

\subsection{Results}

We coded the above procedure in interval arithmetic,
based on Fredrik Johansson's library Arb \cite{arb}, using $53$ bits
of internal precision and $C=10$.
In every case that we tested, we found exactly one consistent choice
of $N$, $w$, and $L_2(A,T)$.

\begin{example}
Consider the genus 2 curve
$$
\href{http://www.lmfdb.org/Genus2Curve/Q/3732/b/477696/1}{y^2+(x^3+x+1)y=-x^6+6x^4+10x^3-33x^2-14x+3}
$$
of minimal discriminant $\Delta_{\min}=477696=2^9\cdot3\cdot311$. Running
through all positive integers~$N$ with $3\cdot311\mid N\mid\Delta_{\min}$, $w=\pm1$,
and all choices of $L_2(A,T)$, we find
that for $N=3732=2^2\cdot3\cdot311$, $w=1$, and $L_2(A,T)=1-T+T^2$,
the value of $S(2^{\frac14}\sqrt{N})-wS(2^{-\frac14}\sqrt{N})$
lies in an interval with center very near $0$ and radius less
than $1.24\times10^{-11}$. (Most of the uncertainty in this
value comes from the error term \eqref{SCerror}.)
For every other choice of $N$, $w$ or $L_2(A,T)$, we find, provably, that
$|S(2^{\frac14}\sqrt{N})-wS(2^{-\frac14}\sqrt{N})|>3.3\times10^{-7}$.
\end{example}

\begin{remark}
One can treat \emph{all} of the $a_n$
as unknowns (not only the ones for even $n$), subject
only to the assumptions that they are integer valued and satisfy the
Ramanujan bound.  This was carried out by Farmer--Koutsoliotas--Lemurell 
\cite{fkl}, who found all solutions to \eqref{LJfunceq}
with $N\le 500$. Assuming the Hasse--Weil conjecture, this includes a
complete classification, up to isogeny, of the abelian surfaces with
conductor $\le500$.  
% (Note that it does not necessarily classify them up to
% isomorphism. Indeed, our search turned up many examples of genus 2 curves with
% conductor less than 500, including some with very large discriminant;
% it is unlikely that we have found all such curves.)
\end{remark}

\begin{remark}
It is possible for a genus 2 curve to have bad reduction at a prime $p$ that does not divide the conductor, so $p \mid \Delta_{\min}$ but $p \nmid N$.  
When this occurs $\Delta_{\min}$ is typically divisible by a large power of $p$ (typically $p^{12}$), and consequently much larger than $10^6$; it is thus not surprising that none of the curves in our database have bad reduction at a prime that does not divide the conductor.
%However, because all of the genus 2 curves in our database have discriminants that are
%either very small,
% or very smooth (divisible by only $2$ and $3$),  
%it is difficult to see this behavior manifest; none of our curves with discriminant $<10^6$ have this property.
% , but there are $55$ curves with $3$-smooth discriminant that have bad reduction at both $2$ and $3$ but whose conductors are either a power of $2$ or a power of $3$.  An example is: 
% \[ y^2 = -2x^5 - 4x^4 + 20x^3 - 8x^2 - x + 2 \]
% with minimal discriminant $2^{21} 3^{12}$ and conductor $2^{10}$.  However, its Jacobian appears to be isogenous to that of the curve $y^2 = -4x^5 + 36x^4 - 72x^3 + 56x^2 - 18x + 2$ with discriminant $2^{31}$.  
% We are led to wonder: if $X$ is a genus $2$ curve with Jacobian $\Jac(X)$ of conductor $N$, is it true that for all primes $p \mid \Delta$ such that $p \nmid N$ there exists a genus $2$ curve $X'$ with good reduction at $p$ such that $\Jac(X')$ is isogenous to $\Jac(X)$?
\end{remark}

\section{Endomorphism data}\label{sec:end}

%In this section, we describe the computation of data about the endomorphism ring of the Jacobian of a genus $2$ curve.  

%\subsection{General strategy}

Let $X:y^2=f(x)$ be a simplified Weierstrass model of a genus 2 curve over~$\Q$, and let $A\colonequals\Jac(X)$ be its Jacobian.  The endomorphism ring $\End(A)$ is an important invariant of the curve (recall from section \ref{sec:notation} that the elements of $\End(A)$ are defined over $\Q$): for example, the presence of nontrivial idempotents in $\End(A)$ indicates that $A$ splits as a product of elliptic curves.  More coarsely, the $\Q$-algebra $\End (A)_{\Q}$ contains information about the isogeny class of~$A$, and the presence of nontrivial idempotents in this algebra indicates whether $A$ is $\Q$-isogenous to a product of elliptic curves. Finally, if $K$ is the minimal field over which all the endomorphisms of $A$ are defined, then as mentioned in section \ref{sec:st_group}, the $\R$-algebra $\End (A_K)_{\R}$ determines and is determined by the identity component $\ST^0 (A)$ of the Sato--Tate group of $A$.  

Our approach to determining $\End (A)$ is to first determine $\End (A_{\Qbar})$, and thereby $\End (A_{\Qbar})_{\Q}$ and $\End (A_{\Qbar})_{\R}$.  We have $\End (A_{\Qbar}) = \End (A_{\C})$,
% \jv{Mumford cite?}\jrs{Why? Since we know that the endomorphism ring over $\C$ is finitely generated as an abelian group, this seems OK.}
so we employ transcendental methods to determine these rings.  (Consequently, our results are not fully rigorous, as they could in principle depend on the precision of our calculations; this is something we hope to rectify in future work.)  From this, we obtain a candidate for the minimal number field $K$ over which all endomorphisms of~$A$ are defined, and for each subfield $F \subset K$ we can then compute $\End(A_F)$ by taking Galois invariants. In other words, we can determine the structure of $\End (A_{\Qbar})$ as a $\Gal (\Qbar / \Q)$-module. 
% Additionally, knowing this structure allows us to determine the Sato--Tate group of $A$ over
% any number field, using the results in (ref), without doing any direct
% calculation over these number fields. We managed in particular to do this for a
% curve whose endomorphism ring is defined over a degree $48$ number field.

\subsection{Reduction to linear algebra}

We apply the techniques of van Wamelen \cite{vW1}.  The Abel--Jacobi map gives a holomorphic isomorphism $A(\C) \simeq V / \Lambda$, where $V = H^0 (X, \omega_X)\spcheck$ is the dual vector space of the global differentials of our genus $g=2$ curve $X$, identified with the tangent space of $A$ at $0$, and $\Lambda \subset V$ is given by $H_1 (X, \Z)$ (via the embedding coming from integration of forms along paths). Concretely, after choosing bases, $\Lambda$ can be embedded as a lattice in $\C^g$ using the \defi{periods} of $X$, that is, the integrals $\int_{\delta_j} \omega_i$, where $\omega_1,\ldots,\omega_g$ is a basis for $H^0 (X, \omega_X)$, and the $\delta_1,\ldots,\delta_{2g}$ is a basis for $H_1 (X, \Z)$. These periods can be evaluated by numerical methods; an especially fast technique when $g=2$ uses the arithmetic-geometric mean and is due to Dupont \cite{Dupont}.

We now find the endomorphisms of this abelian variety using the LLL algorithm \cite{LLL}. To do this, we consider the lattice $\Lambda$ not as a subset of $\C^g$, but rather view it as the standard lattice $\Z^{2g} \subset \R^{2g}$ equipped with a \defi{complex structure}, that is, a map $J : \R^{2g} \to \R^{2g}$ such that $J^2 = -1$. The elements of $\End (A_{\Qbar})$ can then be identified with the maps $R : \Z^{2g} \to \Z^{2g}$ that satisfy
\begin{equation}\label{eq:LLL}
  R J = J R .
\end{equation}
Equation \eqref{eq:LLL} is linear in the integer entries of $R$, and given an approximation of $A$, a basis of approximate solutions of this equation can be found using the LLL algorithm. In practice, we apply this algorithm after determining the period matrix of $X$ (and with it, the matrix $J$) up to $200$ digits of precision. Having in this way found a basis $R_1, \dots , R_d$ of the ring $\End (A_{\Qbar})$ in terms of matrices $R_i : \Z^{2g} \to \Z^{2g}$, we then apply standard meataxe algorithms \cite{Holt} to determine the algebras $\End (A_{\Qbar})_{\Q}$. We may then describe the ring $\End (A_{\Qbar})$ as an order in $\End (A_{\Qbar})_{\Q}$.

As a sanity check, we verify that (to high precision) the structure of the real endomorphism algebra $\End (A_{\Qbar})_{\R}$ agrees with that implied by the Sato--Tate group of $X$, which we (heuristically) computed as described in section \ref{sec:st_group} using entirely different methods.

\subsection{Galois module structure}

The ring $\End (A_{\Qbar})$ is a $\Gal (\Qbar / \Q)$-module in a natural way, and $\End (A)$ can be recovered by taking Galois invariants, that is, by determining the elements of $\End (A_{\Qbar})$ that are fixed under the Galois action. This reduces the question to linear algebra over $\Q$.

%The putative basis of $\End (A_{\Qbar})$ found by using the methods in the previous section comes equipped with an action on the tangent space of $A$ at $0$. Moreover, 
As our curve $X$ is defined over $\Q$, the tangent space of $A_{\Q}$ at $0$ comes equipped with a natural action of $\Gal (\Qbar / \Q)$. The space $H^0 (X_{\Qbar}, \omega_{X_{\Qbar}})$ that is dual to this tangent space has a basis $\mathrm{d}x / y, x\,\mathrm{d}x/y$ defined over $\Q$. We use this to determine the Galois action on $\End (A_{\Qbar})$; denoting the period matrix of $A$ with respect to the given basis of differentials by $\Pi$, the equality
\begin{equation}\label{eq:RtoM}
  M \Pi = \Pi R
\end{equation}
allows us to convert any basis element $R_i : \Z^4 \to \Z^4$ above into a matrix $M_i : \C^2 \to \C^2$. By the argument above, the entries of the matrices $M_i$ lie in $\Qbar$, and $\Gal(\Qbar/\Q)$ acts on $M_i$ entry-wise. Using LLL once more, we can determine a putative Galois number field $L$ over which all these entries are defined. And given a subfield $K$ of $L$, we can determine $\End (A_K)$ by solving a system of linear equations: we find combinations $\sum_{i = 1}^d n_i M_i$ that are stable under all generators $\sigma$ of the corresponding Galois subgroup, that is, rational solutions $n_i$ to the system of equations
\begin{equation}\label{eq:galend}
  \sum_{i = 1}^d n_i (M_i^{\sigma} - M_i) = 0,
\end{equation}
with $\sigma$ running over the aforementioned generators. To determine the endomorphism ring over~$K$, we in fact need a basis of the integral solutions of \eqref{eq:galend}, but this can be efficiently recovered from the basis over $\Q$ using lattice saturation techniques.

\subsection{Decomposition up to isogeny}

The situation where $\End (A_{\Qbar})_{\Q}$ contains a nontrivial idempotent $e$ is of particular interest.  We show in this section how to obtain a splitting of $A \stackrel{\sim}{\rightarrow} E \times E'$ into a product of two elliptic curves, up to isogeny.

The idempotent $e$ gives rise to a matrix $M_e : \C^2 \to \C^2$ whose image has dimension $1$.  In practice, this means that we restrict to a single factor of $\C^2$, one in which the image of the lattice $\Lambda$ under $M_e$ does not become trivial. We thus obtain an isogeny
\begin{equation*}
  L_e : \C^2 / \Lambda \to \C / \Lambda_e,
\end{equation*}
where $\Lambda_e$ is a lattice in $\C$. Using the values $g_4 (\Lambda_e), g_6 (\Lambda_e)$ of the classical elliptic functions, we obtain an elliptic curve $E$ whose complex points are given by $\C / \Lambda_e$. What is more, the properties of these functions together with our $\Q$-rational choice of basis of $H^0 (X_{\C}, \omega_{X_{\C}})$ ensure that the map $L_e$ is defined over a field that is as small as possible. Indeed, in all cases where $e$ was defined over $\Q$ we have obtained rational elliptic curves $E$ whose primes of bad reduction are a subset of those of $X$.

An online implementation of the methods discussed in this section in a mixture of \textsc{Magma} \cite{magma}, \textsc{Pari} \cite{pari} and \textsc{SageMath} \cite{sagemath} can be found at \cite{SijsGitHub}. We conclude by giving some examples.

\subsection{Examples}

\begin{example}
  We first consider the curve
  \begin{center}
    $X\colon$ \href{http://www.lmfdb.org/Genus2Curve/Q/686/a/686/1}{$y^2 + (x^2 + x) y = x^5 + x^4 + 2 x^3 + x^2 + x$}.
  \end{center}
  An approximation of the period matrix for this curve with respect to the basis $\left\{\mathrm{d}x/y, x\,\mathrm{d}x/y \right\}$ is given by
  \begin{equation*}
    \begin{pmatrix}
      3.938 + 1.418 i &
      0.072 + 9.184 i &
      -1.933 - 2.650 i &
      1.933 - 2.650 i \\
      -2.005 - 6.534 i &
      -5.871 - 4.069 i &
      -1.933 + 2.650 i &
      1.933 + 2.650 i
    \end{pmatrix}.
  \end{equation*}
  Using the LLL methods sketched above, we find that with respect to the chosen homology basis, the endomorphism ring is generated as a $\Z$-module by the three elements
  \begin{equation*}
    \begin{pmatrix}
      1 & 0 & 0 & 0 \\
      0 & 1 & 0 & 1 \\
      0 & 0 & 1 & 0 \\
      0 & 0 & 0 & 1
    \end{pmatrix} ,
    \begin{pmatrix}
      0 & 1 & 0 & 1  \\
      1 & 0 & -1 & 0 \\
      0 & 0 & 0 & 1  \\
      0 & 0 & 1 & 0 
    \end{pmatrix} ,
    \begin{pmatrix}
      -1 & 1 & -2 & 3 \\
       2 & -2 & 2 & -4 \\
       1 & -1 & 0 & -1 \\
      -1 & 1 & 0 & 1
    \end{pmatrix} .
  \end{equation*}
  Via equality \eqref{eq:RtoM}, these correspond to the complex matrices
  \begin{equation*}
    \begin{pmatrix}
      1.000 & 0.000 \\
      0.000 & 1.000
    \end{pmatrix} ,
    \begin{pmatrix}
      0.000 & -1.000 \\
      -1.000 & 0.000
    \end{pmatrix} ,
    \begin{pmatrix}
      -0.500 + 1.323 i & -0.500 + 1.323 i \\
      -0.500 + 1.323 i & -0.500 + 1.323 i
    \end{pmatrix},
  \end{equation*}
which are approximations of the following matrices over $\Q (\sqrt{-7})$:
  \begin{equation*}
    \begin{pmatrix}
      1 & 0 \\
      0 & 1
    \end{pmatrix} ,
    \begin{pmatrix}
      0 & -1 \\
      -1 & 0
    \end{pmatrix} ,
    \begin{pmatrix}
      (-1 + \sqrt{-7})/2 & (-1 + \sqrt{-7})/2 \\
      (-1 + \sqrt{-7})/2 & (-1 + \sqrt{-7})/2
    \end{pmatrix} .
  \end{equation*}
The solutions to the equations \eqref{eq:galend} are easy to spot in this case: the elements of $\End (A)$ are exactly those in the span of the first two matrices.

It follows that $\End (A)_{\Q}$ is isomorphic to $\Q \times \Q$; in this algebra, $\End (A)$ corresponds to an order of index $2$ in $\Z \times \Z$. A more careful study shows that $\End (A_{\Qbar})_{\Q} \simeq \Q \times \Q(\sqrt{-7})$, and $\End (A_{\Qbar})$ is an order of index $4$ in $\Z\times \Z[\frac{1+\sqrt{-7}}{2}]$.

We see that the Jacobian $A$ is already split over $\Q$. Constructing idempotents, we obtain $2$-isogenies to the rational elliptic curves with LMFDB labels \href{http://www.lmfdb.org/EllipticCurve/Q/14/a/5}{\texttt{14.a5}} and \href{http://www.lmfdb.org/EllipticCurve/Q/49/a/4}{\texttt{49.a4}}. This agrees with the fact that the conductor of $X$ equals $686 = 2 \cdot 7^3$.
\end{example}

\begin{example}
  Our database contains only two curves with quaternionic multiplication (QM), both of conductor $262144=2^{18}$, with equations \href{http://www.lmfdb.org/Genus2Curve/Q/262144/d/524288/1}{$y^2 = x^5 - x^4 + 4x^3 - 8x^2 + 5x - 1$} and \href{http://www.lmfdb.org/Genus2Curve/Q/262144/d/524288/2}{$y^2 = x^5 + x^4 + 4x^3 + 8x^2 + 5x + 1$}; their Jacobians are isogenous quadratic twists and both admit QM by a maximal order in the rational quaternion algebra of discriminant $6$.
\end{example}

\begin{acknowledgements}\label{ackref}
The authors thank Masha Vlasenko for hosting and the American Institute of Mathematics for funding an LMFDB workshop to facilitate this work, and also thank William Stein for advising us on the use of Google's compute engine and Jan Steffen M\"uller for making us aware of the paper of Berry \cite{Be90}.  Finally, the authors thank the anonymous referees for their comments and corrections.
\end{acknowledgements}

\affiliationone{
   Andrew R. Booker\\
	 School of Mathematics\\
	 University of Bristol\\
	 University Walk\\
	 Bristol, BS8 1TW, UK
   \email{andrew.booker@bristol.ac.uk}}
% Important: Do not put any empty line here.
\affiliationtwo{
   Jeroen Sijsling and John Voight\\
   Department of Mathematics\\
   Dartmouth College\\
   6188 Kemeny Hall\\
   Hanover, NH 03755, USA
   \email{sijsling@gmail.com, jvoight@gmail.com}}
% Important: Do not put any empty line here.
\affiliationthree{
    Andrew V. Sutherland\\
    Department of Mathematics\\
    Massachusetts Institute of Technology\\
    77 Massachusetts Avenue\\
    Cambridge, MA \ 02139, USA
    \email{drew@math.mit.edu}}    
% Important: Do not put any empty line here.
\affiliationfour{%
   Dan Yasaki\\
   Department of Mathematics and Statistics\\
   University of North Carolina at Greensboro\\      
   317 College Avenue\\
   Greensboro, NC 27412, USA
   \email{d\_yasaki@uncg.edu}}

\vfill\newpage
\oneappendix % use \appendix if you have more than one appendix
\section{Tables}

%We conclude with some tables giving some information about our data set.  

\begin{center}
\vspace{-1ex}
\begin{table}[h!]
\setlength{\tabcolsep}{8pt}
\begin{tabular}{llrrrr}
&&\multicolumn{4}{c}{discriminant bound}\\
\cmidrule{3-6}
size & box & $10^3$ & $10^4$ & $10^5$ & $10^6$\\
\midrule
$\approx 10^{11}$ & $S_1(12)$ & \num{47} & \num{921} & \num{8301} & \num{56724}\\
$\approx 10^{12}$ & $S_1(17)$ & \num{50} & \num{953} & \num{8622} & \num{59482}\\
$\approx 10^{13}$ & $S_1(24)$ & \num{52} & \num{979} & \num{8852} & \num{61401}\\
$\approx 10^{14}$ & $S_1(33)$ & \num{53} & \num{987} & \num{8993} & \num{62696}\\
$\approx 10^{15}$ & $S_1(46)$ & \num{54} & \num{1007} & \num{9121} & \num{63697}\\
$\approx 10^{16}$ & $S_1(64)$ & \num{54} & \num{1020} & \num{9232} & \num{64524}\\
$\approx 10^{17}$ & $S_1(90)$ & \num{54} & \num{1028} & \num{9304} & \num{65093}\\
\midrule
$\approx 10^{11}$ & $S_2(2,1.83)$ & \num{50} & \num{944} & \num{8486} & \num{58127}\\
$\approx 10^{12}$ & $S_2(2,2.02)$ & \num{51} & \num{967} & \num{8761} & \num{60546}\\
$\approx 10^{13}$ & $S_2(2,2.27)$ & \num{52} & \num{986} & \num{8921} & \num{62061}\\
$\approx 10^{14}$ & $S_2(2,2.53)$ & \num{52} & \num{1000} & \num{9053} & \num{63150}\\
$\approx 10^{15}$ & $S_2(2,2.83)$ & \num{54} & \num{1015} & \num{9150} & \num{64035}\\
$\approx 10^{16}$ & $S_2(2,3.16)$ & \num{54} & \num{1023} & \num{9217} & \num{64537}\\
$\approx 10^{17}$ & $S_2(2,3.51)$ & \num{56} & \num{1027} & \num{9274} & \num{64923}\\
\midrule
$\approx 10^{11}$ & $S_3(3.00)$ & \num{48} & \num{937} & \num{8404} & \num{57402}\\
$\approx 10^{12}$ & $S_3(3.47)$ & \num{49} & \num{958} & \num{8590} & \num{59043}\\
$\approx 10^{13}$ & $S_3(4.00)$ & \num{50} & \num{978} & \num{8801} & \num{60859}\\
$\approx 10^{14}$ & $S_3(4.69)$ & \num{50} & \num{982} & \num{8870} & \num{61644}\\
$\approx 10^{15}$ & $S_3(5.37)$ & \num{52} & \num{999} & \num{9003} & \num{62740}\\
$\approx 10^{16}$ & $S_3(6.17)$ & \num{54} & \num{1007} & \num{9106} & \num{63543}\\
$\approx 10^{17}$ & $S_3(7.14)$ & \num{54} & \num{1010} & \num{9167} & \num{64021}\\
\bottomrule
\end{tabular}
\caption{Searching in boxes for genus 2 curves $y^2+h(x)y=f(x)$ with $|\Delta(f,h)|\le 10^n$.}\label{tab:boxsearch}
\end{table}
\end{center}

\begin{center}
\begin{table}[h!]
\begin{tabular}{lrrrr}
set & cardinality & models & curves & discs \\\midrule
$S_1$ & $1.02\times 10^{17}$ & \num{1908760} & \num{65093} &\num{56132}\\\midrule
$S_2$ & $9.84\times 10^{16}$ & \num{1295050} &\num{64923} & \num{56009}\\
$S_1\cup S_2$ & $2.00\times 10^{17}$ & \num{2147485} &\num{65638} & \num{56549}\\\midrule
$S_3$ & $1.01\times 10^{17}$ & \num{1299092} & \num{64021} & \num{55240}\\
$S_1\cup S_2\cup S_3$ & $3.01\times 10^{17}$ & \num{2218870} &\num{65892} & \num{56728}\\\midrule
$S_4$ &$2.10\times 10^{16}$ & \num{1453127} & \num{64768} & \num{55820}\\
$S_1\cup S_2\cup S_3\cup S_4$ & $3.22\times 10^{17}$ &\num{2232548} & \num{66130} & \num{56907}\\\bottomrule
\end{tabular}
\caption{Searching for integral models with discriminants $|\Delta(f,h)|\le 10^6$.}\label{tab:boxsearchfinal}
% in the\\ ${}\,\,\ \ \qquad\qquad$sets $S_1 \colonequals S_1(90)$, $S_2 \colonequals S_2(2,3.51)$, $S_3 \colonequals S_3(7.14)$, $S_4 \colonequals S_4(10,10)$.}\label{tab:boxsearchfinal}
\end{table}
\end{center}
%
%\begin{center}
%\begin{table}[h!]
%\begin{tabular}{lll}
%LMFDB label & minimal model & Sato-Tate group\\\midrule
%\href{http://www.lmfdb.org/Genus2Curve/Q/4096/b/65536/1}{4096.b.65536.1} & $y^2=x^5-x$ & $J(C_2)$\\
%\href{http://www.lmfdb.org/Genus2Curve/Q/5184/a/46656/1}{5184.a.46656.1} & $y^2+x^3y=x^3+2$ & $J(C_2)$\\
%\href{http://www.lmfdb.org/Genus2Curve/Q/40000/e/200000/1}{40000.e.200000.1} & $y^2+x^3y=x^5-5x^3-10x^2-8x-2$ & $J(C_4)$\\
%\href{http://www.lmfdb.org/Genus2Curve/Q/20736/i/373248/1}{20736.i.373248.1} & $y^2+x^3y=-2$ & $D_{2,1}$\\
%\href{http://www.lmfdb.org/Genus2Curve/Q/20736/k/373248/1}{20736.k.373248.1} & $y^2+x^3y=2$ & $D_{2,1}$\\
%\href{http://www.lmfdb.org/Genus2Curve/Q/65536/a/65536/1}{65536.a.65536.1} & $y^2=x^5+x$ & $D_{2,1}$\\
%\href{http://www.lmfdb.org/Genus2Curve/Q/2916/b/11664/1}{2916.b.11664.1} & $y^2+y=x^6$ & $D_{3,2}$\\
%\href{http://www.lmfdb.org/Genus2Curve/Q/11664/a/11664/1}{11664.a.11664.1} & $y^2+y=-x^6$ & $D_{6,2}$\\
%\bottomrule
%\end{tabular}
%\caption{Some interesting curves. \avs{What distinguishes these?}}  
%\end{table}
%\end{center}

\begin{center}
\begin{table}[h!]
\begin{tabular}{rrrrrr}
$\#W(\Q)$ & total & $\GL_2$-type & $\Q$-simple & $\Qbar$-simple\\ \midrule
0 & \num{32616} & \num{2265} & \num{30462}& \num{30326} \\
1 & \num{24611} & 46 &  \num{24580} & \num{24574} \\
2 & \num{8005} & 490 & \num{7522} & \num{7461} \\
3 & \num{886} & 35 & 881 & 851\\
4 & \num{40} & 10 & 32 & 20\\
\toprule
 & \num{66158} & \num{2846} & \num{63477} &\num{63232} \\ 
\bottomrule
\end{tabular}
\caption{Number of $\Q$-rational Weierstrass points $\#W(\Q)$.}
\end{table}
\end{center}

\vfill\newpage

\begin{center}
\begin{table}[h!]
\begin{tabular}{rlrrrr}
$\#G$ & $G$ & total & $\GL_2$-type & $\Q$-simple & $\Qbar$-simple\\ \midrule
1 & $\{0\}$ & \href{http://www.lmfdb.org/Genus2Curve/Q/?&torsion=%5B%5D}{\num{44190}} & \href{http://www.lmfdb.org/Genus2Curve/Q/?&torsion=%5B%5D&is_gl2_type=True}{586} & \num{43694}& \num{43664} \\
2 & $\Z/2\Z$ & \href{http://www.lmfdb.org/Genus2Curve/Q/?&torsion=%5B2%5D}{\num{14681}} & \href{http://www.lmfdb.org/Genus2Curve/Q/?&torsion=%5B2%5D&is_gl2_type=True}{848} & \num{13845}&\num{13737} \\
3 & $\Z/3\Z$ & \href{http://www.lmfdb.org/Genus2Curve/Q/?&torsion=%5B3%5D}{\num{2295}} & \href{http://www.lmfdb.org/Genus2Curve/Q/?&torsion=%5B3%5D&is_gl2_type=True}{279} & \num{2019} & \num{2006} \\
4 & $\Z/2\Z\times \Z/2\Z$ & \href{http://www.lmfdb.org/Genus2Curve/Q/?&torsion=%5B2%2C2%5D}{\num{1352}} & \href{http://www.lmfdb.org/Genus2Curve/Q/?&torsion=%5B2%2C2%5D&is_gl2_type=True}{244} & \num{1143} & \num{1085}  \\
4 & $\Z/4\Z$ & \href{http://www.lmfdb.org/Genus2Curve/Q/?&torsion=%5B4%5D}{\num{1402}} & \href{http://www.lmfdb.org/Genus2Curve/Q/?&torsion=%5B4%5D&is_gl2_type=True}{338} & \num{1065} & \num{1063} \\
5 & $\Z/5\Z$ & \href{http://www.lmfdb.org/Genus2Curve/Q/?&torsion=%5B5%5D}{725} & \href{http://www.lmfdb.org/Genus2Curve/Q/?&torsion=%5B5%5D&is_gl2_type=True}{84} & 650 & 650 \\
6 & $\Z/6\Z$ & \href{http://www.lmfdb.org/Genus2Curve/Q/?&torsion=%5B6%5D}{595} & \href{http://www.lmfdb.org/Genus2Curve/Q/?&torsion=%5B6%5D&is_gl2_type=True}{210} & 387 & 378 \\
7 & $\Z/7\Z$ & \href{http://www.lmfdb.org/Genus2Curve/Q/?&torsion=%5B7%5D}{97} & \href{http://www.lmfdb.org/Genus2Curve/Q/?&torsion=%5B7%5D&is_gl2_type=True}{13} & 86 & 86\\
8 & $\Z/8\Z$ & \href{http://www.lmfdb.org/Genus2Curve/Q/?&torsion=%5B8%5D}{201} & \href{http://www.lmfdb.org/Genus2Curve/Q/?&torsion=%5B8%5D&is_gl2_type=True}{28} & 175 & 174 \\
8 & $\Z/2\Z\times\Z/4\Z$ & \href{http://www.lmfdb.org/Genus2Curve/Q/?&torsion=%5B2%2C4%5D}{159} & \href{http://www.lmfdb.org/Genus2Curve/Q/?&torsion=%5B2%2C4%5D&is_gl2_type=True}{92} & 67 & 66 \\
8 & $\Z/2\Z\times\Z/2\Z\times\Z/2\Z$ & \href{http://www.lmfdb.org/Genus2Curve/Q/?&torsion=%5B2%2C2%2C2%5D}{33} & \href{http://www.lmfdb.org/Genus2Curve/Q/?&torsion=%5B2%2C2%2C2%5D&is_gl2_type=True}{5} & 30 & 18\\
9 & $\Z/9\Z$ & \href{http://www.lmfdb.org/Genus2Curve/Q/?&torsion=%5B9%5D}{30} & \href{http://www.lmfdb.org/Genus2Curve/Q/?&torsion=%5B9%5D&is_gl2_type=True}{1} & 30 & 29\\
9 & $\Z/3\Z\times\Z/3\Z$ & \href{http://www.lmfdb.org/Genus2Curve/Q/?&torsion=%5B3%2C3%5D}{8} & \href{http://www.lmfdb.org/Genus2Curve/Q/?&torsion=%5B3%2C3%5D&is_gl2_type=True}{5} & 2 & 0\\
10 & $\Z/10\Z$ & \href{http://www.lmfdb.org/Genus2Curve/Q/?&torsion=%5B10%5D}{131} & \href{http://www.lmfdb.org/Genus2Curve/Q/?&torsion=%5B10%5D&is_gl2_type=True}{5} & 127 & 126\\
11 & $\Z/11\Z$ & \href{http://www.lmfdb.org/Genus2Curve/Q/?&torsion=%5B11%5D}{8} & \href{http://www.lmfdb.org/Genus2Curve/Q/?&torsion=%5B11%5D&is_gl2_type=True}{2} & 8 & 8\\
12 & $\Z/12\Z$ & \href{http://www.lmfdb.org/Genus2Curve/Q/?&torsion=%5B12%5D}{59} & \href{http://www.lmfdb.org/Genus2Curve/Q/?&torsion=%5B12%5D&is_gl2_type=True}{40} & 19 & 18\\
12 & $\Z/2\Z\times\Z/6\Z$ & \href{http://www.lmfdb.org/Genus2Curve/Q/?&torsion=%5B2%2C6%5D}{66} & \href{http://www.lmfdb.org/Genus2Curve/Q/?&torsion=%5B2%2C6%5D&is_gl2_type=True}{22} & 47 & 44\\
13 & $\Z/13\Z$ & \href{http://www.lmfdb.org/Genus2Curve/Q/?&torsion=%5B13%5D}{7} & \href{http://www.lmfdb.org/Genus2Curve/Q/?&torsion=%5B13%5D&is_gl2_type=True}{0} & 7 & 7\\
14 & $\Z/14\Z$ & \href{http://www.lmfdb.org/Genus2Curve/Q/?&torsion=%5B14%5D}{12} & \href{http://www.lmfdb.org/Genus2Curve/Q/?&torsion=%5B14%5D&is_gl2_type=True}{1} & 12 & 12\\
15 & $\Z/15\Z$ & \href{http://www.lmfdb.org/Genus2Curve/Q/?&torsion=%5B15%5D}{17} & \href{http://www.lmfdb.org/Genus2Curve/Q/?&torsion=%5B15%5D&is_gl2_type=True}{7} & 10 & 10\\
16 & $\Z/16\Z$ & \href{http://www.lmfdb.org/Genus2Curve/Q/?&torsion=%5B16%5D}{3} & \href{http://www.lmfdb.org/Genus2Curve/Q/?&torsion=%5B16%5D&is_gl2_type=True}{0} & 3 & 3\\
16 & $\Z/2\Z\times\Z/8\Z$ & \href{http://www.lmfdb.org/Genus2Curve/Q/?&torsion=%5B2%2C8%5D}{31} & \href{http://www.lmfdb.org/Genus2Curve/Q/?&torsion=%5B2%2C8%5D&is_gl2_type=True}{11} & 20 & 20\\
16 & $\Z/4\Z\times \Z/4\Z$ & \href{http://www.lmfdb.org/Genus2Curve/Q/?&torsion=%5B4%2C4%5D}{3} & \href{http://www.lmfdb.org/Genus2Curve/Q/?&torsion=%5B4%2C4%5D&is_gl2_type=True}{2} & 0 & 0\\
16 & $\Z/2\Z\times\Z/2\Z\times\Z/4\Z$ & \href{http://www.lmfdb.org/Genus2Curve/Q/?&torsion=%5B2%2C2%2C4%5D}{3} & \href{http://www.lmfdb.org/Genus2Curve/Q/?&torsion=%5B2%2C2%2C4%5D&is_gl2_type=True}{2} & 1 & 1\\
17 & $\Z/17\Z$ & \href{http://www.lmfdb.org/Genus2Curve/Q/?&torsion=%5B17%5D}{1} & \href{http://www.lmfdb.org/Genus2Curve/Q/?&torsion=%5B17%5D&is_gl2_type=True}{0} & 1 & 1\\
18 & $\Z/18\Z$ & \href{http://www.lmfdb.org/Genus2Curve/Q/?&torsion=%5B18%5D}{3} & \href{http://www.lmfdb.org/Genus2Curve/Q/?&torsion=%5B18%5D&is_gl2_type=True}{0} & 3 & 3\\
18 & $\Z/3\Z\times\Z/6\Z$ & \href{http://www.lmfdb.org/Genus2Curve/Q/?&torsion=%5B3%2C6%5D}{6} & \href{http://www.lmfdb.org/Genus2Curve/Q/?&torsion=%5B3%2C6%5D&is_gl2_type=True}{5} & 0 & 0\\
19 & $\Z/19\Z$ & \href{http://www.lmfdb.org/Genus2Curve/Q/?&torsion=%5B19%5D}{1} & \href{http://www.lmfdb.org/Genus2Curve/Q/?&torsion=%5B19%5D&is_gl2_type=True}{1} & 1 & 0\\
20 & $\Z/20\Z$ & \href{http://www.lmfdb.org/Genus2Curve/Q/?&torsion=%5B20%5D}{6} & \href{http://www.lmfdb.org/Genus2Curve/Q/?&torsion=%5B20%5D&is_gl2_type=True}{1} & 5 & 5\\
20 & $\Z/2\Z\times\Z/10\Z$ & \href{http://www.lmfdb.org/Genus2Curve/Q/?&torsion=%5B2%2C10%5D}{7} & \href{http://www.lmfdb.org/Genus2Curve/Q/?&torsion=%5B2%2C10%5D&is_gl2_type=True}{2} & 6 & 5\\
21 & $\Z/21\Z$ & \href{http://www.lmfdb.org/Genus2Curve/Q/?&torsion=%5B21%5D}{5} & \href{http://www.lmfdb.org/Genus2Curve/Q/?&torsion=%5B21%5D&is_gl2_type=True}{3} & 3 & 2\\
22 & $\Z/22\Z$ & \href{http://www.lmfdb.org/Genus2Curve/Q/?&torsion=%5B22%5D}{2} & \href{http://www.lmfdb.org/Genus2Curve/Q/?&torsion=%5B22%5D&is_gl2_type=True}{0} & 2 & 2\\
24 & $\Z/24\Z$ & \href{http://www.lmfdb.org/Genus2Curve/Q/?&torsion=%5B24%5D}{4} & \href{http://www.lmfdb.org/Genus2Curve/Q/?&torsion=%5B24%5D&is_gl2_type=True}{3} & 1 & 1\\
24 & $\Z/2\Z\times\Z/12\Z$ & \href{http://www.lmfdb.org/Genus2Curve/Q/?&torsion=%5B2%2C12%5D}{4} & \href{http://www.lmfdb.org/Genus2Curve/Q/?&torsion=%5B2%2C12%5D&is_gl2_type=True}{2} & 2 & 2\\
24 & $\Z/2\Z\times\Z/2\Z\times\Z/6\Z$ & \href{http://www.lmfdb.org/Genus2Curve/Q/?&torsion=%5B2%2C2%2C6%5D}{2} & \href{http://www.lmfdb.org/Genus2Curve/Q/?&torsion=%5B2%2C2%2C6%5D&is_gl2_type=True}{1} & 1 & 1\\
27 & $\Z/27\Z$ & \href{http://www.lmfdb.org/Genus2Curve/Q/?&torsion=%5B27%5D}{1} & \href{http://www.lmfdb.org/Genus2Curve/Q/?&torsion=%5B27%5D&is_gl2_type=True}{0} & 1 & 1\\
27 & $\Z/3\Z\times\Z/9\Z$ & \href{http://www.lmfdb.org/Genus2Curve/Q/?&torsion=%5B3%2C9%5D}{1} & \href{http://www.lmfdb.org/Genus2Curve/Q/?&torsion=%5B3%2C9%5D&is_gl2_type=True}{1} & 0 & 0\\
28 & $\Z/28\Z$ & \href{http://www.lmfdb.org/Genus2Curve/Q/?&torsion=%5B28%5D}{1} & \href{http://www.lmfdb.org/Genus2Curve/Q/?&torsion=%5B28%5D&is_gl2_type=True}{0} & 1 & 1\\
28 & $\Z/2\Z\times\Z/14\Z$ & \href{http://www.lmfdb.org/Genus2Curve/Q/?&torsion=%5B2%2C14%5D}{1} & \href{http://www.lmfdb.org/Genus2Curve/Q/?&torsion=%5B2%2C14%5D&is_gl2_type=True}{0} & 1 & 1\\
29 & $\Z/29\Z$ & \href{http://www.lmfdb.org/Genus2Curve/Q/?&torsion=%5B29%5D}{1} & \href{http://www.lmfdb.org/Genus2Curve/Q/?&torsion=%5B29%5D&is_gl2_type=True}{0} & 1 & 1\\
32 & $\Z/2\Z\times\Z/2\Z\times\Z/8\Z$ & \href{http://www.lmfdb.org/Genus2Curve/Q/?&torsion=%5B2%2C2%2C8%5D}{2} & \href{http://www.lmfdb.org/Genus2Curve/Q/?&torsion=%5B2%2C2%2C8%5D&is_gl2_type=True}{2} & 0 & 0\\
36 & $\Z/6\Z\times\Z/6\Z$ & \href{http://www.lmfdb.org/Genus2Curve/Q/?&torsion=%5B6%2C6%5D}{1} & \href{http://www.lmfdb.org/Genus2Curve/Q/?&torsion=%5B6%2C6%5D&is_gl2_type=True}{0} & 0 & 0\\
39 & $\Z/39\Z$ & \href{http://www.lmfdb.org/Genus2Curve/Q/?&torsion=%5B39%5D}{1} &  \href{http://www.lmfdb.org/Genus2Curve/Q/?&torsion=%5B39%5D&is_gl2_type=True}{0} & 1 & 1\\
\midrule
 &  & \num{66158} & \num{2846} &\num{63477} & \num{63232} \\
\bottomrule
\end{tabular}
\caption{Torsion subgroups.}
\end{table}
\end{center}

\vfill\newpage

\begin{center}
\begin{table}[h!]
\begin{tabular}{llrrrrrr}
$\Aut(X)$ & $\Aut(X_{\Qbar})$ & total & $\GL_2$-type & $\Q$-simple & $\Qbar$-simple\\ \midrule
$C_2$ & $C_2$ & \href{http://www.lmfdb.org/Genus2Curve/Q/?aut_grp=%5B2%2C+1%5D&geom_aut_grp=%5B2%2C+1%5D}{\num{63310}} & \href{http://www.lmfdb.org/Genus2Curve/Q/?aut_grp=%5B2%2C+1%5D&geom_aut_grp=%5B2%2C+1%5D&is_gl2_type=True}{174} & \num{63234} & \num{63227} \\
$C_2$ & $V_4$ & \href{http://www.lmfdb.org/Genus2Curve/Q/?aut_grp=%5B2%2C+1%5D&geom_aut_grp=%5B4%2C+2%5D}{125} & \href{http://www.lmfdb.org/Genus2Curve/Q/?aut_grp=%5B2%2C+1%5D&geom_aut_grp=%5B4%2C+2%5D&is_gl2_type=True}{1} & 125 & 0\\
$C_2$ & $D_8$ & \href{http://www.lmfdb.org/Genus2Curve/Q/?aut_grp=%5B2%2C+1%5D&geom_aut_grp=%5B8%2C+3%5D}{17} & \href{http://www.lmfdb.org/Genus2Curve/Q/?aut_grp=%5B2%2C+1%5D&geom_aut_grp=%5B8%2C+3%5D&is_gl2_type=True}{0} & 17 & 0\\
$C_2$ & $C_{10}$ & \href{http://www.lmfdb.org/Genus2Curve/Q/?aut_grp=%5B2%2C+1%5D&geom_aut_grp=%5B10%2C+2%5D}{5} & \href{http://www.lmfdb.org/Genus2Curve/Q/?aut_grp=%5B2%2C+1%5D&geom_aut_grp=%5B10%2C+2%5D&is_gl2_type=True}{0} & 5 & 5\\
$C_2$ & $D_{12}$ & \href{http://www.lmfdb.org/Genus2Curve/Q/?aut_grp=%5B2%2C+1%5D&geom_aut_grp=%5B12%2C+4%5D}{23} & \href{http://www.lmfdb.org/Genus2Curve/Q/?aut_grp=%5B2%2C+1%5D&geom_aut_grp=%5B12%2C+4%5D&is_gl2_type=True}{1} & 23 & 0\\
\midrule
$C_4$ & $D_8$ & \href{http://www.lmfdb.org/Genus2Curve/Q/?aut_grp=%5B4%2C+1%5D&geom_aut_grp=%5B8%2C+3%5D}{13} & \href{http://www.lmfdb.org/Genus2Curve/Q/?aut_grp=%5B4%2C+1%5D&geom_aut_grp=%5B8%2C+3%5D&is_gl2_type=True}{13} & 13 & 0\\
$C_4$ & $\tilde{S}_4$ & \href{http://www.lmfdb.org/Genus2Curve/Q/?aut_grp=%5B4%2C+1%5D&geom_aut_grp=%5B48%2C+29%5D}{2} & \href{http://www.lmfdb.org/Genus2Curve/Q/?aut_grp=%5B4%2C+1%5D&geom_aut_grp=%5B48%2C+29%5D&is_gl2_type=True}{2} & 2 & 0\\
\midrule
$V_4$ & $V_4$ & \href{http://www.lmfdb.org/Genus2Curve/Q/?aut_grp=%5B4%2C+2%5D&geom_aut_grp=%5B4%2C+2%5D}{\num{2573}} & \href{http://www.lmfdb.org/Genus2Curve/Q/?aut_grp=%5B4%2C+2%5D&geom_aut_grp=%5B4%2C+2%5D&is_gl2_type=True}{\num{2573}} & 0 & 0\\
$V_4$ & $D_8$ & \href{http://www.lmfdb.org/Genus2Curve/Q/?aut_grp=%5B4%2C+2%5D&geom_aut_grp=%5B8%2C+3%5D}{16} & \href{http://www.lmfdb.org/Genus2Curve/Q/?aut_grp=%5B4%2C+2%5D&geom_aut_grp=%5B8%2C+3%5D&is_gl2_type=True}{16} & 0 & 0\\
$V_4$ & $D_{12}$ & \href{http://www.lmfdb.org/Genus2Curve/Q/?aut_grp=%5B4%2C+2%5D&geom_aut_grp=%5B12%2C+4%5D}{3} & \href{http://www.lmfdb.org/Genus2Curve/Q/?aut_grp=%5B4%2C+2%5D&geom_aut_grp=%5B12%2C+4%5D&is_gl2_type=True}{3} & 0 & 0\\
$V_4$ & $2D_{12}$ & \href{http://www.lmfdb.org/Genus2Curve/Q/?aut_grp=%5B4%2C+2%5D&geom_aut_grp=%5B24%2C+8%5D}{4} & \href{http://www.lmfdb.org/Genus2Curve/Q/?aut_grp=%5B4%2C+2%5D&geom_aut_grp=%5B24%2C+8%5D&is_gl2_type=True}{4} & 0 & 0\\
$V_4$ & $\tilde{S}_4$ & \href{http://www.lmfdb.org/Genus2Curve/Q/?aut_grp=%5B4%2C+2%5D&geom_aut_grp=%5B48%2C+29%5D}{1} & \href{http://www.lmfdb.org/Genus2Curve/Q/?aut_grp=%5B4%2C+2%5D&geom_aut_grp=%5B48%2C+29%5D&is_gl2_type=True}{1} & 0 & 0\\
\midrule
$C_6$ & $D_{12}$ & \href{http://www.lmfdb.org/Genus2Curve/Q/?aut_grp=%5B6%2C+2%5D&geom_aut_grp=%5B12%2C+4%5D}{58} & \href{http://www.lmfdb.org/Genus2Curve/Q/?aut_grp=%5B6%2C+2%5D&geom_aut_grp=%5B12%2C+4%5D&is_gl2_type=True}{58} & 58 & 0\\
\midrule
$D_8$ & $D_8$ & \href{http://www.lmfdb.org/Genus2Curve/Q/?aut_grp=%5B8%2C+3%5D&geom_aut_grp=%5B8%2C+3%5D}{4} & \href{http://www.lmfdb.org/Genus2Curve/Q/?aut_grp=%5B8%2C+3%5D&geom_aut_grp=%5B8%2C+3%5D&is_gl2_type=True}{0} & 0 & 0\\
\midrule
$D_{12}$ & $D_{12}$ & \href{http://www.lmfdb.org/Genus2Curve/Q/?aut_grp=%5B12%2C+4%5D&geom_aut_grp=%5B12%2C+4%5D}{4} & \href{http://www.lmfdb.org/Genus2Curve/Q/?aut_grp=%5B12%2C+4%5D&geom_aut_grp=%5B12%2C+4%5D&is_gl2_type=True}{0} & 0 & 0\\
\midrule
&  & \num{66158} & \num{2846} &\num{63477} & \num{63232} \\
\bottomrule
\end{tabular}
\caption{Automorphism groups.}
\end{table}
\end{center}

\vfill\newpage

\begin{center}
\begin{table}[h!]
\begin{tabular}{llrrrr}
type & $\ST$ & total & $\GL_2$-type & $\Q$-simple & $\Qbar$-simple\\ \midrule
\textbf{A} & $\USp(4)$ & \href{http://www.lmfdb.org/Genus2Curve/Q/?st_group=USp%284%29}{\num{63107}} & \href{http://www.lmfdb.org/Genus2Curve/Q/?st_group=USp%284%29&is_gl2_type=True}{0} & \num{63107} & \num{63107}\\
\midrule
\textbf{B} & $G_{3,3}$& \href{http://www.lmfdb.org/Genus2Curve/Q/?st_group=G_%7B3%2C3%7D}{\num{2440}} &  \href{http://www.lmfdb.org/Genus2Curve/Q/?st_group=G_%7B3%2C3%7D&is_gl2_type=True}{\num{2440}} & \num{97} & \num{97} \\
\textbf{B} & $N(G_{3,3})$ & \href{http://www.lmfdb.org/Genus2Curve/Q/?st_group=N%28G_%7B3%2C3%7D%29}{144} & \href{http://www.lmfdb.org/Genus2Curve/Q/?st_group=N%28G_%7B3%2C3%7D%29&is_gl2_type=True}{0} & 144 & 19\\
\midrule
\textbf{C} & $N(G_{1,3})$ & \href{http://www.lmfdb.org/Genus2Curve/Q/?st_group=N%28G_%7B1%2C3%7D%29}{303} & \href{http://www.lmfdb.org/Genus2Curve/Q/?st_group=N%28G_%7B1%2C3%7D%29&is_gl2_type=True}{303} & 0 & 0 \\
\midrule
\textbf{D} & $F_{ac}$ & \href{http://www.lmfdb.org/Genus2Curve/Q/?st_group=F_%7Bac%7D}{6} & \href{http://www.lmfdb.org/Genus2Curve/Q/?st_group=F_%7Bac%7D&is_gl2_type=True}{0} & 6 & 6\\
\textbf{D} & $F_{a,b}$ & \href{http://www.lmfdb.org/Genus2Curve/Q/?st_group=F_%7Ba%2Cb%7D}{0} & \href{http://www.lmfdb.org/Genus2Curve/Q/?st_group=F_%7Ba%2Cb%7D&is_gl2_type=True}{0} & 0 & 0\\
\midrule
\textbf{E} & $E_1$ & \href{http://www.lmfdb.org/Genus2Curve/Q/?st_group=E_1}{8} & \href{http://www.lmfdb.org/Genus2Curve/Q/?st_group=E_1&is_gl2_type=True}{0} & 0 & 0\\
\textbf{E} & $E_2$ & \href{http://www.lmfdb.org/Genus2Curve/Q/?st_group=E_2}{3} & \href{http://www.lmfdb.org/Genus2Curve/Q/?st_group=E_2&is_gl2_type=True}{3} & 3 & 0\\
\textbf{E} & $E_3$ & \href{http://www.lmfdb.org/Genus2Curve/Q/?st_group=E_3}{7} & \href{http://www.lmfdb.org/Genus2Curve/Q/?st_group=E_3&is_gl2_type=True}{7} & 7 & 0\\
\textbf{E} & $E_4$ & \href{http://www.lmfdb.org/Genus2Curve/Q/?st_group=E_4}{10} & \href{http://www.lmfdb.org/Genus2Curve/Q/?st_group=E_4&is_gl2_type=True}{10} & 10 & 0\\
\textbf{E} & $E_6$ & \href{http://www.lmfdb.org/Genus2Curve/Q/?st_group=E_6}{51} & \href{http://www.lmfdb.org/Genus2Curve/Q/?st_group=E_6&is_gl2_type=True}{51} & 51 & 0\\
\textbf{E} & $J(E_1)$ & \href{http://www.lmfdb.org/Genus2Curve/Q/?st_group=J%28E_1%29}{24} & \href{http://www.lmfdb.org/Genus2Curve/Q/?st_group=J%28E_1%29&is_gl2_type=True}{24} & 2 & 0\\
\textbf{E} & $J(E_2)$ & \href{http://www.lmfdb.org/Genus2Curve/Q/?st_group=J%28E_2%29}{9} & \href{http://www.lmfdb.org/Genus2Curve/Q/?st_group=J%28E_2%29&is_gl2_type=True}{0} & 9 & 2\\
\textbf{E} & $J(E_3)$ & \href{http://www.lmfdb.org/Genus2Curve/Q/?st_group=J%28E_3%29}{4} & \href{http://www.lmfdb.org/Genus2Curve/Q/?st_group=J%28E_3%29&is_gl2_type=True}{0} & 4 & 0\\
\textbf{E} & $J(E_4)$ & \href{http://www.lmfdb.org/Genus2Curve/Q/?st_group=J%28E_4%29}{17} & \href{http://www.lmfdb.org/Genus2Curve/Q/?st_group=J%28E_4%29&is_gl2_type=True}{0} & 17 & 1\\
\textbf{E} & $J(E_6)$ & \href{http://www.lmfdb.org/Genus2Curve/Q/?st_group=J%28E_6%29}{17} & \href{http://www.lmfdb.org/Genus2Curve/Q/?st_group=J%28E_6%29&is_gl2_type=True}{0} & 17 & 0\\
\midrule
\textbf{F} & $J(C_2)$ & \href{http://www.lmfdb.org/Genus2Curve/Q/?st_group=J%28C_2%29}{2} & \href{http://www.lmfdb.org/Genus2Curve/Q/?st_group=J%28C_2%29&is_gl2_type=True}{2} & 2 & 0\\
\textbf{F} & $J(C_4)$ & \href{http://www.lmfdb.org/Genus2Curve/Q/?st_group=J%28C_4%29}{1} & \href{http://www.lmfdb.org/Genus2Curve/Q/?st_group=J%28C_4%29&is_gl2_type=True}{1} & 0 & 0\\
\textbf{F} & $J(C_6)$ & \href{http://www.lmfdb.org/Genus2Curve/Q/?st_group=J%28C_6%29}{0} & \href{http://www.lmfdb.org/Genus2Curve/Q/?st_group=J%28C_6%29&is_gl2_type=True}{0} & 0 & 0\\
\textbf{F} & $J(D_2)$ & \href{http://www.lmfdb.org/Genus2Curve/Q/?st_group=J%28D_2%29}{0} & \href{http://www.lmfdb.org/Genus2Curve/Q/?st_group=J%28D_2%29&is_gl2_type=True}{0} & 0 & 0\\
\textbf{F} & $J(D_3)$ & \href{http://www.lmfdb.org/Genus2Curve/Q/?st_group=J%28D_3%29}{0} & \href{http://www.lmfdb.org/Genus2Curve/Q/?st_group=J%28D_3%29&is_gl2_type=True}{0} & 0 & 0\\
\textbf{F} & $J(D_4)$ & \href{http://www.lmfdb.org/Genus2Curve/Q/?st_group=J%28D_4%29}{0} & \href{http://www.lmfdb.org/Genus2Curve/Q/?st_group=J%28D_4%29&is_gl2_type=True}{0} & 0 & 0\\
\textbf{F} & $J(D_6)$ & \href{http://www.lmfdb.org/Genus2Curve/Q/?st_group=J%28D_6%29}{0} & \href{http://www.lmfdb.org/Genus2Curve/Q/?st_group=J%28D_6%29&is_gl2_type=True}{0} & 0 & 0\\
\textbf{F} & $J(T)$ & \href{http://www.lmfdb.org/Genus2Curve/Q/?st_group=J%28T%29}{0} & \href{http://www.lmfdb.org/Genus2Curve/Q/?st_group=J%28T%29&is_gl2_type=True}{0} & 0 & 0\\
\textbf{F} & $J(O)$ & \href{http://www.lmfdb.org/Genus2Curve/Q/?st_group=J%28O%29}{0} & \href{http://www.lmfdb.org/Genus2Curve/Q/?st_group=J%28O%29&is_gl2_type=True}{0} & 0 & 0\\
\textbf{F} & $C_{2,1}$ & \href{http://www.lmfdb.org/Genus2Curve/Q/?st_group=C_%7B2%2C1%7D}{0} & \href{http://www.lmfdb.org/Genus2Curve/Q/?st_group=C_%7B2%2C1%7D&is_gl2_type=True}{0} & 0 & 0\\
\textbf{F} & $C_{6,1}$ & \href{http://www.lmfdb.org/Genus2Curve/Q/?st_group=C_%7B6%2C1%7D}{0} & \href{http://www.lmfdb.org/Genus2Curve/Q/?st_group=C_%7B6%2C1%7D&is_gl2_type=True}{0} & 0 & 0\\
\textbf{F} & $D_{2,1}$ & \href{http://www.lmfdb.org/Genus2Curve/Q/?st_group=D_%7B2%2C1%7D}{3} & \href{http://www.lmfdb.org/Genus2Curve/Q/?st_group=D_%7B2%2C1%7D&is_gl2_type=True}{3} & 0 & 0\\
\textbf{F} & $D_{4,1}$ & \href{http://www.lmfdb.org/Genus2Curve/Q/?st_group=D_%7B4%2C1%7D}{0} & \href{http://www.lmfdb.org/Genus2Curve/Q/?st_group=D_%7B4%2C1%7D&is_gl2_type=True}{0} & 0 & 0\\
\textbf{F} & $D_{6,1}$ & \href{http://www.lmfdb.org/Genus2Curve/Q/?st_group=D_%7B6%2C1%7D}{0} & \href{http://www.lmfdb.org/Genus2Curve/Q/?st_group=D_%7B6%2C1%7D&is_gl2_type=True}{0} & 0 & 0\\
\textbf{F} & $D_{3,2}$ & \href{http://www.lmfdb.org/Genus2Curve/Q/?st_group=D_%7B3%2C2%7D}{1} & \href{http://www.lmfdb.org/Genus2Curve/Q/?st_group=D_%7B3%2C2%7D&is_gl2_type=True}{1} & 0 & 0\\
\textbf{F} & $D_{4,2}$ & \href{http://www.lmfdb.org/Genus2Curve/Q/?st_group=D_%7B4%2C2%7D}{0} & \href{http://www.lmfdb.org/Genus2Curve/Q/?st_group=D_%7B4%2C2%7D&is_gl2_type=True}{0} & 0 & 0\\
\textbf{F} & $D_{6,2}$ & \href{http://www.lmfdb.org/Genus2Curve/Q/?st_group=D_%7B6%2C2%7D}{1} & \href{http://www.lmfdb.org/Genus2Curve/Q/?st_group=D_%7B6%2C2%7D&is_gl2_type=True}{1} & 0 & 0\\
\textbf{F} & $O_1$ & \href{http://www.lmfdb.org/Genus2Curve/Q/?st_group=O_1}{0} & \href{http://www.lmfdb.org/Genus2Curve/Q/?st_group=O_1&is_gl2_type=True}{0} & 0 & 0\\
\midrule
& & \num{66158} & \num{2846} &\num{63477} & \num{63232} \\
\bottomrule
\end{tabular}
\caption{Sato-Tate groups.}
\end{table}
\end{center}

\vfill\newpage

\begin{center}
\begin{table}[h!]
\begin{tabular}{rrrrr}
$2$-Selmer rank & $\Jac(X)[2]$ rank & analytic rank & $\Sha[2]$ rank & total\\
\midrule
0 & 0 & 0 & 0 & \num{3878}\\
1 & 0 & 0 & 1 & \num{604}\\
1 & 0 & 1 & 0 & \num{20191}\\
1 & 1 & 0 & 0 & \num{4701}\\
\midrule
2 & 0 & 0 & 2 & \num{445}\\
2 & 0 & 1 & 1 & \num{680}\\
2 & 0 & 2 & 0 & \num{18481}\\
2 & 1 & 0 & 1 & \num{477}\\
2 & 1 & 1 & 0 & \num{8434}\\
2 & 2 & 0 & 0 & \num{886}\\
\midrule
3 & 0 & 0 & 3 & \num{4}\\
3 & 0 & 1 & 2 & \num{173}\\
3 & 0 & 2 & 1 & \num{88}\\
3 & 0 & 3 & 0 & \num{2832}\\
3 & 1 & 0 & 2 & \num{908}\\
3 & 1 & 1 & 1 & \num{307}\\
3 & 1 & 2 & 0 & \num{1958}\\
3 & 2 & 0 & 1 & \num{44}\\
3 & 2 & 1 & 0 & \num{529}\\
3 & 3 & 0 & 0 & \num{36}\\
\midrule
4 & 0 & 2 & 2 & \num{2}\\
4 & 0 & 4 & 0 & \num{10}\\
4 & 1 & 0 & 3 & \num{24}\\
4 & 1 & 1 & 2 & \num{217}\\
4 & 1 & 2 & 1 & \num{16}\\
4 & 1 & 3 & 0 & \num{44}\\
4 & 2 & 0 & 2 & \num{101}\\
4 & 2 & 1 & 1 & \num{17}\\
4 & 2 & 2 & 0 & \num{13}\\
4 & 3 & 1 & 0 & \num{4}\\
\midrule
5 & 1 & 0 & 4 & \num{9}\\
5 & 1 & 1 & 3 & \num{7}\\
5 & 1 & 2 & 2 & \num{1}\\
5 & 1 & 3 & 1 & \num{1}\\
5 & 2 & 0 & 3 & \num{12}\\
5 & 2 & 1 & 2 & \num{17}\\
5 & 2 & 2 & 1 & \num{2}\\
\midrule
6 & 1 & 1 & 4 & \num{2}\\
6 & 2 & 0 & 4 & \num{2}\\
6 & 2 & 1 & 3 & \num{1}\\
\midrule
1.617 & 0.309 & 1.215 & 0.093 & \num{66158} \\
\bottomrule
\end{tabular}
\caption{$2$-Selmer rank, torsion rank, analytic rank, $\Sha$ $2$-rank, and averages.}
\end{table}
\end{center}

\begin{remark}
The fact that the average analytic rank in our database is significantly greater than the expected average value of $0.5$ is likely due to the fact that having small coefficients made it more likely that $4f+h^2$ will be a sextic whose leading coefficient is a square, in which case the difference of the two points at infinity is a rational divisor that typically has infinite order (it is torsion if and only if the continued fraction expansion of $y$ in the hyperelliptic function field $\Q[x,y[/(y^2-4f-h^2)$ is periodic \cite{Be90}). When this occurs the Mordell-Weil rank of the Jacobian is $1$ greater than it would be generically.  
\end{remark}


\begin{thebibliography}{99}% Replace 9 by 99 if 10 or more references
%
% Please note the use of "\and" between author names below
%
% \bibitem{Goddard}
% {\bibname P. Goddard, A. Kent \and D. I. Olive}, `Unitary
% representations of the Virasoro and Supervirasoro algebras', {\em
% Comm. Math. Phys. }103 (1986) 105.

\newcommand{\mytextit}[1]{{`#1'}}
\newcommand{\mytextbf}[1]{#1}

\bibitem{BGGP05}
{\bibname M. H. Baker, E. Gonz\'alez-Jim\'enez, J. Gonz\'alez, \and B. Poonen}, \href{http://muse.jhu.edu/journals/american_journal_of_mathematics/v127/127.6baker.pdf}{\mytextit{Finiteness results for modular curves of genus at least $2$}}, {\em Amer. J. Math.} \mytextbf{127} (2005), 1325--1387.
%\bibitem{BK12}
%{\bibname G. Banaszak \and K. S. Kedlaya}, \href{http://www.iumj.indiana.edu/IUMJ/FULLTEXT/2015/64/5438}{\mytextit{An algebraic Sato--Tate group and Sato--Tate conjecture}}, {\em Indiana Univ. Math. J.} \mytextbf{64} (2015), 245--274.
\bibitem{Be90}
{\bibname T. G. Berry},
\href{http://link.springer.com/article/10.1007/BF01191166}{\mytextit{On periodicity of continued fractions in hyperelliptic function fields}}, Arch. Math. \textbf{55} (1990), 259--266.
\bibitem{AntwerpIV}
{\bibname B. J. Birch \and W. Kuyk (eds.)}, \href{http://link.springer.com/book/10.1007/BFb0097580}{{\em Modular functions of one variable IV}}, Lecture Notes in Math. \mytextbf{476}, Springer-Verlag, Berlin, 1975.
\bibitem{B05}
{\bibname A. R. Booker}, \href{http://www.maths.bris.ac.uk/~maarb/public/papers/modularity.pdf}{\mytextit{Numerical tests of modularity}}, {\em J. Ramanujan Math. Soc.} \mytextbf{20} (2005), No. 4, 283--339.
\bibitem{MagmaHandbook}
{\bibname W. Bosma, J. Cannon, C. Fieker, \and A. Steel (eds.)}, \href{http://magma.maths.usyd.edu.au/magma/handbook/text/1457}{{\em Handbook of Magma functions}}, Edition 2.16 (2010), 5017 pages. 
\bibitem{magma}
{\bibname W. Bosma, J. Cannon, \and C. Playoust}, \href{http://www.sciencedirect.com/science/article/pii/S074771719690125X}{\mytextit{The Magma algebra system I: The user language}},  {\em J. Symbolic Comput.} \mytextbf{24} (1997), 235--265.
\bibitem{BK94}
{\bibname A. Brumer \and K. Kramer}, \href{http://www.numdam.org/item?id=CM_1994__92_2_227_0}{\mytextit{The conductor of an abelian variety}}, {\em Compos. Math.} \mytextbf{92} (1994), 227--248.
\bibitem{BK14}
{\bibname A. Brumer \and K. Kramer}, \href{http://www.ams.org/journals/tran/2014-366-05/S0002-9947-2013-05909-0}{\mytextit{Paramodular abelian varieties of odd conductor}}, {\em Trans. Amer. Math. Soc.} \mytextbf{366} (2014), 2463--2516.
\bibitem{BouwWewers}
{\bibname I. I. Bouw \and S. Wewers}, \href{http://arxiv.org/abs/1211.4459}{\mytextit{Computing $L$-functions and semistable reduction of superelliptic curves}}, arXiv:1211.4459v4, 2014.
\bibitem{CQNP2}
{\bibname G. Cardona, E. Nart, \and J. Pujolas}, \href{http://link.springer.com/article/10.1007/s00209-004-0750-0}{\mytextit{Curves of genus two over fields of even characteristic}}, {\em Math. Z.} \mytextbf{250} (2005), 177--201.
\bibitem{CQNP1}
{\bibname G. Cardona \and J. Quer}, \href{http://www.worldscientific.com/doi/abs/10.1142/9789812701640_0006}{\mytextit{Field of moduli and field of definition for curves of genus 2}}, in {\em Computational aspects of algebraic curves}, Lecture Notes Ser. Comput. \mytextbf{13}, World Sci. Publ. (2005), Hackensack, NJ, 71--83.
\bibitem{Cremona0}
{\bibname J. E. Cremona}, \href{http://johncremona.github.io/ecdata/}{{\em Algorithms for modular elliptic curves}}, 2nd ed., Cambridge Univ. Press, Cambridge, 1997. 
\bibitem{CremonaTables}
{\bibname J. E. Cremona}, \href{http://link.springer.com/chapter/10.1007\%2F11792086_2}{\mytextit{The elliptic curve database for conductors to 130000}}, in {\em Algorithmic Number Theory (ANTS VII)} (F. Hess, S. Pauli, M. Pohst eds.), LNCS \mytextbf{4076}, Springer, 2006, 11--29.
\bibitem{Cr14}
{\bibname J. E. Cremona}, \href{https://homepages.warwick.ac.uk/staff/J.E.Cremona/ftp/data/}{\mytextit{Elliptic curve data}}, at \url{https://homepages.warwick.ac.uk/staff/J.E.Cremona/ftp/data/}.
\bibitem{D04}
{\bibname T. Dokchitser}, \href{http://projecteuclid.org/euclid.em/1090350929}{\mytextit{Computing special values of motivic $L$-functions}}, {\em Exp. Math.} \mytextbf{13} (2004), 137--149.
\bibitem{Dupont}
{\bibname R. Dupont}, \href{http://www.lix.polytechnique.fr/Labo/Regis.Dupont/these_soutenance.pdf}{\mytextit{Moyenne arithm\'etico-g\'eom\'etrique, suites de Borchardt et applications}}, Ph.D. thesis, Laboratoire d'Informatique de l'\'Ecole Polytechnique, at \url{http://www.lix.polytechnique.fr/Labo/Regis.Dupont/these_soutenance.pdf} 2006.
\bibitem{fkl}
{\bibname D. Farmer, S. Koutsoliotas, \and S. Lemurell}, \href{http://arxiv.org/abs/1502.00850}{\mytextit{Varieties via their $L$-functions}}, arXiv:1502.00850v1, 2015.
\bibitem{FKRS12}
{\bibname F. Fit\'e, K. S. Kedlaya, V. Rotger, \and A. V. Sutherland}, \href{http://dx.doi.org/10.1112/S0010437X12000279}{\mytextit{Sato--Tate distributions and Galois endomorphism modules in genus $2$}}, {\em Compos. Math.} \mytextbf{148} (2012) 1390--1442.
%\bibitem{FS14}
%{\bibname F. Fit\'e \and A. V. Sutherland}, \href{http://msp.org/ant/2014/8-3/p02.xhtml}{\mytextit{Sato--Tate distributions of twists of $y^2=x^5-x$ and $y^2=x^6+1$}}, {\em Algebra Number Theory} \mytextbf{8} (2014), 543--585.
\bibitem{GGR05}
{\bibname J. Gonz\'alez, J. Gu\`ardia, \and V. Rotger}, \href{http://journals.impan.pl/cgi-bin/doi?aa116-3-3}{\mytextit{Abelian surfaces of $\GL_2$-type as Jacobians of curves}}, {\em Acta Arith.} \mytextbf{116} (2005), 263--287.
\bibitem{GR04}
{\bibname J. Gonz\'alez \and V. Rotger}, \href{http://dx.doi.org/10.1155/S1073792804131826}{\mytextit{Equations of Shimura curves of genus two}}, {\em Int. Math. Res. Not.} \mytextbf{14} (2004), 661--674.
\bibitem{GG02}
{\bibname E. Gonz\'alez-Jim\'enez \and J. Gonz\'alez}, \href{http://www.ams.org/journals/mcom/2003-72-241/S0025-5718-02-01458-8}{\mytextit{Modular curves of genus $2$}}, {\em Math. Comp.} \mytextbf{72} (2002), 397--418.
\bibitem{GGG02}
{\bibname E. Gonz\'alez-Jim\'enez, J. Gonz\'alez, \and J. Gu\`ardia}, \href{http://link.springer.com/chapter/10.1007\%2F3-540-45455-1_15}{\mytextit{Computations on modular Jacobian surfaces}}, in {\em Algorithmic Number Theory (ANTS V)} (C. Fieker, D.R. Kohel eds.), LNCS \mytextbf{2369}, Springer (2002), 189--197.
\bibitem{Google}
Google Inc., \href{https://cloud.google.com/}{\mytextit{Google Cloud Platform}}, at \url{https://cloud.google.com}.
\bibitem{HS14}
{\bibname D. Harvey \and A. V. Sutherland}, \href{http://dx.doi.org/10.1112/S1461157014000187}{\mytextit{Computing Hasse-Witt matrices of hyperelliptic curves in average polynomial time}}, in {\em Algorithmic Number Theory (ANTS XI)}, LMS J. Comput. Math. \mytextbf{17} (2014), 257--273.
\bibitem{HS16}
{\bibname D. Harvey \and A. V. Sutherland}, \href{http://arxiv.org/abs/1410.5222}{\mytextit{Computing Hasse-Witt matrices of hyperelliptic curves in average polynomial time II}}, in {\em Frobenius Distributions: Lang-Trotter and Sato-Tate Conjectures} (D. Kohel and I. E. Shparlinski eds.), Contemporary Mathematics \textbf{663}, AMS, 127--148, to appear.
\bibitem{Holt}
{\bibname D. Holt}, \href{http://ebooks.cambridge.org/chapter.jsf?bid=CBO9780511565830&cid=CBO9780511565830A016}{\mytextit{The Meataxe as a tool in computational group theory}}, in {\em The atlas of finite groups: ten years on (Birmingham, 1995)}, London Math. Soc. Lecture Note Ser. \mytextbf{249}, Cambridge Univ. Press (1998), 74--81.
\bibitem{Igusa}
{\bibname J.-I. Igusa}, \href{http://www.jstor.org/stable/1970233}{\mytextit{Arithmetic variety of moduli for genus two}}, {\em Ann. Math.} \mytextbf{72} (1960), 612--649.
\bibitem{arb}
{\bibname F. Johansson}, \href{http://fredrikj.net/arb/}{\mytextit{Arb: a C library for ball arithmetic}}, {\em ACM Communications in Computer Algebra} \mytextbf{47} (2013), No. 4, 166--169.
\bibitem{KS08}
{\bibname K. S. Kedlaya \and A. V. Sutherland}, \href{http://link.springer.com/chapter/10.1007/978-3-540-79456-1_21}{\mytextit{Computing $L$-series of hyperelliptic curves}}, in {\em Algorithmic Number Theory (ANTS VIII)} (A. J. van der Poorten, A. Stein eds.), LNCS \mytextbf{5011} Springer (2008), Berlin, 312--326.
\bibitem{LLL}
{\bibname A. K. Lenstra, H. W. Lenstra, Jr. \and L. Lov\' asz}, \href{http://link.springer.com/article/10.1007\%2FBF01457454}{\mytextit{Factoring polynomials with rational coefficients}}, {\em Math. Ann.} \mytextbf{261} (1982), 513--534.
\bibitem{Li94}
{\bibname Q. Liu}, \href{http://www.numdam.org/item?id=CM_1994__94_1_51_0}{\mytextit{Conducteur et discriminant minimal de courbes de genre 2}}, {\em Compos. Math.} \mytextbf{94} (1994), 51--79.
\bibitem{Li94-alg}
{\bibname Q. Liu}, \href{http://gdz.sub.uni-goettingen.de/dms/load/img/?PPN=GDZPPN00221184X&IDDOC=239088}{\mytextit{Mod\`eles minimaux des courbes de genre deux}}, {\em J. Reine Angew. Math.} \mytextbf{453} (1994), 137--164.
\bibitem{Li96}
{\bibname Q. Liu}, \href{http://www.ams.org/journals/tran/1996-348-11/S0002-9947-96-01684-4/}{\mytextit{Mod\`eles entiers des courbes hyperelliptiques sur un corps valuation discr\`ete}}, {\em Trans. Amer. Math. Society} \mytextbf{348} (1996), 4577--4610.
\bibitem{lmfdb}
The LMFDB Collaboration, \href{http://www.lmfdb.org}{\mytextit{The $L$-functions and modular forms database}}, at \url{http://www.lmfdb.org}.
\bibitem{Lo94}
P. Lockhart, \href{http://www.ams.org/journals/tran/1994-342-02/S0002-9947-1994-1195511-X/}{\mytextit{On the discriminant of a hyperelliptic curve}}, {\em Trans. Amer. Math. Soc.} \mytextbf{342} (1994), 729--752.
\bibitem{MS93}
{\bibname J. R. Merriman \and N. P. Smart}, \href{http://dx.doi.org/10.1017/S030500410007153X}{\mytextit{Curves of genus $2$ with good reduction away from $2$ with a rational Weierstrass point}}, {\em Math. Proc. Cambridge Philos. Soc.} \mytextbf{114} (1993), 203--214.
\bibitem{Mestre}
{\bibname J.-F. Mestre}, \href{http://link.springer.com/chapter/10.1007/978-1-4612-0441-1_21}{\mytextit{Construction de courbes de genre 2 \`a partir de leurs modules}}, in {\em Effective methods in algebraic geometry}, Progr. Math. 94, Birkh\"auser (1991), Boston, 313--334.
\bibitem{pari}
The PARI-Group, \href{http://pari.math.u-bordeaux.fr/}{\mytextit{PARI/GP version {\tt 2.7.0}}}, Bordeaux (2014), available at \url{http://pari.math.u-bordeaux.fr/}.
\bibitem{Po94}
{\bibname B. Poonen}, \href{http://link.springer.com/chapter/10.1007\%2F3-540-61581-4_63}{\mytextit{Computational aspects of curves of genus at least $2$}}, in {\em Algorithmic Number Theory (ANTS II)} (H. Cohen ed.), LNCS \mytextbf{1122}, Springer (1994), Berlin, 283--306.
\bibitem{sagemath}
The Sage Developers, \href{http://www.sagemath.org/}{\mytextit{SageMath, the Sage Mathematics Software System (Version 7.2)}} (2016), available at \url{http://www.sagemath.org/}.
%\bibitem{PY15}
%{\bibname C. Poor \and D. S. Yuen}, \href{http://www.ams.org/journals/mcom/2015-84-293/S0025-5718-2014-02870-6/}{\mytextit{Paramodular cusp forms}}, {\em Math. Comp.} \mytextbf{84} (2015), 1401--1438.
%\bibitem{Si92}
%{\bibname A. Silverberg}, \href{http://www.sciencedirect.com/science/article/pii/0022404992901412}{\mytextit{Fields of definition for homomorphisms of abelian varieties}}, {\em J. Pure Appl. Algebra} \mytextbf{77} (1992), 253--262.
\bibitem{SijsGitHub}
{\bibname J. Sijsling}, \href{https://github.com/JRSijsling/heuristic_endomorphisms/}{\mytextit{Heuristic determination of endomorphism rings of curves}} (2016), available at \url{https://github.com/JRSijsling/heuristic_endomorphisms/}.
\bibitem{smalljac}
{\bibname A. V. Sutherland}, \texttt{smalljac}, version 4.1.3, available at \url{http://math.mit.edu/~drew}.
\bibitem{Sm97}
{\bibname N. P. Smart}, \href{http://dx.doi.org/10.1112/S002461159700035X}{\mytextit{$S$-unit equations, binary forms, and curves of genus $2$}}, {\em Proc. London Math. Soc.} \mytextbf{75} (1997), 271--307.
\bibitem{So04}
{\bibname K. Soundararajan}, \href{http://dx.doi.org/10.4153/CMB-2004-046-0}{\mytextit{Strong multiplicity one for the Selberg class}}, {\em Canad. Math. Bull.} \mytextbf{47} (2004), 468--474.
\bibitem{SteinWatkins}
{\bibname W. Stein \and M. Watkins}, \href{http://link.springer.com/chapter/10.1007\%2F3-540-45455-1_22}{\mytextit{A database of elliptic curves---first report}}, in {\em Algorithmic Number Theory (ANTS V)} (C. Fieker, D. Kohel eds.), LNCS \mytextbf{2369}, Springer (2002), 267--275.
\bibitem{St13}
{\bibname M. Stoll}, \href{http://www.mathe2.uni-bayreuth.de/stoll/data/}{\mytextit{Genus $2$ curve data files}}, at \url{http://www.mathe2.uni-bayreuth.de/stoll/data/}.
%\bibitem{UN73}
%{\bibname Y. Namikawa \and K. Ueno}, \href{https://eudml.org/doc/154162}{\mytextit{The complete classification of fibres in pencils of curves of genus two}}, {\em Manuscripta Math.} \mytextbf{9} (1973), 143--186.
\bibitem{vW1}
{\bibname P. B. van Wamelen}, \href{http://link.springer.com/chapter/10.1007/978-3-540-37634-7_5}{\mytextit{Computing with the analytic Jacobian of a genus 2 curve}}, in {\em Discovering Mathematics with Magma}, Algorithms Comput. Math. \mytextbf{19}, Springer (2006), Berlin, 117--135.

\end{thebibliography}
\end{document}